\documentclass[aop,MSNbibl,seceqn,dvips]{arximspdf}
\usepackage{mathbh}
\usepackage{graphicx}

\doi{10.1214/13-AOP862} %kopijuoti is PTS
\volume{41}
\issue{6}
\pubyear{2013}
\firstpage{4214}
\lastpage{4247}

\makeatletter
\newcommand{\rrvert}{\vert}
\newcommand{\llvert}{\vert}
\newcommand{\eqref}[1]{(\ref{#1})}

\newcommand{\E}{\mathbb{E}}
\newcommand{\R}{\mathbb{R}}
\newcommand{\N}{\mathbb{N}}
\newcommand{\GOE}{{\operatorname{GOE}}}
\newcommand{\Crt}{\mathrm{Crt}}
\def\Ai{\operatorname{Ai} }
\def\d{\mathrm{d}}
\def\<{\langle}
\def\>{\rangle}
\newcommand{\Pro}{\mathbb{P}}
\newcommand{\indi}{\mathbf{1}}
\renewcommand\Re{\operatorname{Re}}
\renewcommand\Im{\operatorname{Im}}
\newcommand{\Var}{\operatorname{Var}}
\newtheorem{theorem}{Theorem}[section]
\newtheorem{question}{Question}[section]
\newtheorem{corollary}{Corollary}[section]
\newtheorem{lemma}{Lemma}%[section]
\newtheorem{proposition}{Proposition}%[section]

\newproclaim{definition}{Definition}[section]
\newproclaim{remark}{Remark}%[section]
\newproclaim{example}{Example}%[section]
\newproclaim{conjecture}{Conjecture}
\makeatother

\begin{document}
\begin{frontmatter}

\title{Complexity of random smooth functions on the high-dimensional sphere}
\runtitle{Complexity of random smooth functions}

\begin{aug}
\author[a]{\fnms{Antonio} \snm{Auffinger}\corref{}\thanksref{t2}\ead[label=e1]{auffing@math.uchicago.edu}} \and
\author[b]{\fnms{Gerard} \snm{Ben Arous}\thanksref{t3}\ead[label=e2]{benarous@cims.nyu.edu}}
\thankstext{t2}{Supported in part by NSF Grants DMS-08-06180 and DMS-05-00923.}
\thankstext{t3}{Supported in part by NSF Grants DMS-08-06180 and OISE-0730136.}
\runauthor{A. Auffinger and G. Ben Arous}
\affiliation{University of Chicago and New York University}
\address[a]{University of Chicago\\
5734 S. University Avenue\\
Chicago, Illinois 60637\\
USA \\
\printead{e1}}

\address[b]{Courant Institute of Mathematical Sciences\\
New York University\\
251 Mercer Street\\
New York, New York 10012\\
USA\\
\printead{e2}}
\end{aug}

% HISTORY:
\received{\smonth{10} \syear{2011}}
\revised{\smonth{3} \syear{2013}}

% ABSTRACT
%
\begin{abstract}
We analyze the landscape of general smooth Gaussian functions on the
sphere in dimension $N$, when $N$ is large. We give an explicit formula
for the asymptotic complexity of the mean number of critical points of
finite and diverging index at any level of energy and for the mean
Euler characteristic of level sets. We then find two possible scenarios
for the bottom landscape, one that has a layered structure of critical
values and a strong correlation between indexes and critical values and
another where even at levels below the limiting ground state energy the
mean number of local minima is exponentially large. We end the paper by
discussing how these results can be interpreted in the language of spin
glasses models.
\end{abstract}

% KEYWORDS
% Pirmas kwd is didziosios raides
%


\end{frontmatter}

%s1 #&#
\section{Introduction}

This work deals with the number of critical points of Gaussian smooth
functions on the $N$ dimensional sphere. The questions addressed in
this paper can be phrased as: What does a random Morse function look
like on a high-dimensional sphere? How many critical values of given
index, or below a given level? What can be said about the topology of
its level sets? We investigate the number of critical points of given
index in level sets below a given value, as well as the topology of the
level sets through their mean Euler characteristic. Our main result is
that these functions have an exponentially large number of critical
points of given index, and that the Euler characteristic of the level
sets have a very interesting oscillatory behavior. Moreover we find an
invariant to distinguish between two very different classes of these
functions that we describe below.

Let us know describe the functions that we will analyze. For $N \geq
1$, let $S^{N-1}(\sqrt N)\subset\mathbb R^N$ be the Euclidean sphere
of radius
$\sqrt N$,
\[
S^{N-1}(\sqrt N): = \Biggl\{ \bolds\sigma=(\sigma_1,
\ldots,\sigma_N) \in\R^N \dvtx\frac{1}{N}\sum
_{i=1}^N \sigma_i^2
= 1 \Biggr\}.
\]

Consider the Gaussian function defined on
$S^{N-1}(\sqrt N)$ by
%
%
%e1.1 #&#
\begin{equation}
\label{Hamiltonian} H_{N,p}(\bolds\sigma) = \frac{1}{N^{(p-1)/2}}
\sum_{i_1, \ldots, i_p=1}^N J_{i_1, \ldots, i_p}
\sigma_{i_1}\cdots\sigma_{i_p},
\end{equation}
where $J_{i_1, \ldots, i_p}$ are independent centered standard Gaussian random
variables.

Equivalently, $H_{N,p}$ is the centered Gaussian process on the
sphere\break
$S^{N-1}(\sqrt N)$ whose covariance is given by
%
%
%e1.2 #&#
\begin{equation}
\mathbb E \bigl[H_{N,p}(\bolds\sigma)H_{N,p}\bigl(
\bolds\sigma'\bigr) \bigr] =N^{1-p} \Biggl(\sum
_{i=1}^N \sigma_i
\sigma'_i \Biggr)^p = N R\bigl(\bolds
\sigma, \bolds\sigma'\bigr)^p,
\end{equation}
where $R$ is the normalized inner product
$R(\bolds\sigma, \bolds\sigma'):=
\frac{1}{N}\< \bolds\sigma, \bolds\sigma'\>=\frac
{1}{N}\sum_{i=1}^N \sigma_i\sigma'_i$.

Given a sequence $\bolds\beta= (\beta_p)_{p\in\N, p\geq2}$ of
positive real numbers such that
%
%
%e1.3 #&#
\begin{equation}
\label{e1b} \sum_{p=2}^{\infty} 2^p
\beta_p < \infty,
\end{equation}
let
%
%
%e1.4 #&#
\begin{equation}
H_{N}(\bolds\sigma) = \sum_{p=2}^{\infty}
\beta_p H_{N,p}(\bolds\sigma),
\end{equation}
where for any pair of values $p\neq p'$, the Hamiltonians $H_{N,p},
H_{N,p'}$ are independent. Condition \eqref{e1b} is more than enough to
guarantee that the above sum is a.s. finite, and the Hamiltonian
$H_{N}$ is a.s. smooth and Morse; see Theorem 11.3.1 of~\cite{AT07}.

In this case, we have that
%
%
%e1.5 #&#
\begin{equation}
\label{duali} \mathbb E \bigl[H_{N}(\bolds\sigma)H_{N}
\bigl(\bolds\sigma'\bigr) \bigr]= N \sum
_{p=2}^\infty\beta_p^2 \bigl(R
\bigl(\bolds\sigma, \bolds\sigma'\bigr)
\bigr)^p = N \nu\bigl(R\bigl(\bolds\sigma, \bolds
\sigma'\bigr)\bigr),
\end{equation}
where
%
%
%e1.6 #&#
\begin{equation}
\label{linsum} \nu(t):= \sum_{p=2}^{\infty}
\beta_p^2 t^p.
\end{equation}

We will fix the variance of $H_N$ by assuming
\[
\nu(1)= \sum_{p=2}^{\infty}
\beta_p^2 = 1.
\]

A word of comment is needed here. By Schoenberg's theorem \cite
{Schoenberg}, if $\nu(R(\bolds\sigma, \bolds\sigma'))$
is a
positive-definite function for all $N$ and all $ \bolds\sigma,
\bolds\sigma' \in S^{N-1}(\sqrt N)$, then $\nu$ can be written as
a linear sum as in \eqref{linsum}. This remark implies that we are
exhausting all possible covariances given as \eqref{duali} that satisfy
\eqref{e1b}. The importance of \eqref{e1b} is to ensure smoothness of
the process $H_N$.

From now on, we call the function $\nu$ a mixture. If $\nu= \beta
_p^2t^p$, for some $p\geq2$, we call $\nu$ a pure mixture.
Note that $\nu$ is smooth with
%
%
%e1.7 #&#
\begin{equation}
\nu'(1):=\nu' \neq0, \qquad \nu''(1)
:= \nu''>0.
\end{equation}

If we consider the random variable $X$ that assigns probability $\beta
_p^2$ to the integer~$p$, then its probability measure is given by $\mu
_X= \sum\beta_p^2 \delta_p$ and
%
%
%e1.8 #&#
\begin{equation}
\label{eq2p} \E X = \nu'\quad \mbox{and}\quad \alpha^2:= \Var X
= \nu'' + \nu' - \nu'^2.
\end{equation}
A mixture is pure if and only if $\alpha=0$. Furthermore, note that
$\nu
''\geq\nu'$ with equality only in the pure case with $p=2$. The
parameters $\nu', \nu''$ and $\alpha^2$ will be fundamental in our analysis.

We now introduce the main object of our study. For any
open set $B\subset\R$ and any integer $0\le k < N$, we consider the
(random) number
$\Crt_{N,k}(B)$ of critical values of the
function $H_{N }$ in the set $NB=\{Nx\dvtx x\in B \}$ with index equal to
$k$,
%
%
%e1.9 #&#
\begin{equation}
\label{defWk} \Crt_{N,k}(B) = \sum_{\bolds\sigma: \nabla
H_{N}(\bolds\sigma) = 0 }
\indi\bigl\{ H_{N}(\bolds\sigma) \in NB\bigr\} \indi\bigl\{ i
\bigl(\nabla^2 H_{N}(\bolds\sigma)\bigr) = k\bigr\}.
\end{equation}

Here $\nabla$, $\nabla^2$ are the gradient and the Hessian restricted to
$S^{N-1}(\sqrt N)$, and $i(\nabla^2 H_{N}(\bolds\sigma))$ is the
number of negative eigenvalues
of the Hessian $\nabla^2 H_{N}$, called the index of the Hessian at
$\bolds\sigma$. We will also consider the total number
$\Crt_{N}(B)$ of critical values of the function $H_{N}$ in the set $NB$
(whatever their index)
%
%
%e1.10 #&#
\begin{equation}
\label{edefW} \Crt_{N}(B) = \sum_{\bolds\sigma: \nabla
H_{N}(\bolds\sigma) = 0 }
\indi\bigl\{ H_{N}(\bolds\sigma) \in NB\bigr\}.
\end{equation}

Our first results will give exact and asymptotic formulas for the mean values
$ \E\Crt_{N,k}(B)$ and $\E\Crt_{N}(B)$, when $N\to\infty$
and $k$, $B$ and $\nu$ are fixed. This initial computation uses the
method developed in \cite{ABC}, where this study was initiated for pure
mixtures.

%
%th1.1 #&#
\begin{theorem}\label{maintheorem} For any fixed integer $k \geq0$,
there exists a continuous function $\theta_{k,\nu}(u)$, called the
$k$-complexity function, explicitly given in \eqref
{complexityfunction2}, such that, for any open set $B \subseteq\R$,
\label{critical2}
%
%
%e1.11 #&#
\begin{equation}
\lim_{N\to\infty} \frac{1}{N} \log\E\Crt_{N,k}(B) =
\sup_{u \in B} \theta_{k,\nu}(u).
\end{equation}
\end{theorem}

%
%f1 #&#
\begin{figure}

\includegraphics{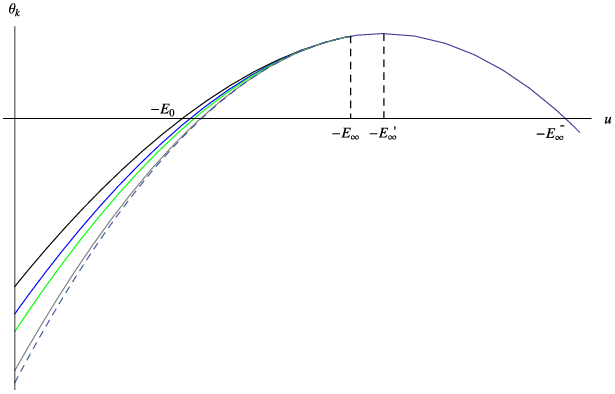}

\caption{$k$-complexity functions $\theta_{k,\nu}(u)$ for $-6\leq u
\leq-1$, $k=1,2,3,5$ in the case where $\nu$ is pure-like,
that is, $\theta_{k,\nu}(-E_{\infty})>0$. The dashed line is the
continuation of the parabola that describes $\theta_{k,\nu}(u)$ in
the interval $[-E_{\infty},\infty)$ where they all agree.}\label{fig1}
\end{figure}

We decide to postpone to Section~\ref{sectionlo} the explicit
expression of the $k$-complexity functions $\theta_{k,\nu}(u)$.
However, we describe some important properties of these functions (see
Figure~\ref{fig1}) in the proposition below. We first fix four important
thresholds that depend on $\nu$. Let
%
%
%e1.12 #&#
\begin{equation}
E'_\infty:= \frac{2 \nu' \sqrt{\nu''}}{\nu'+\nu''}, \qquad E_\infty:=
\frac{2\nu''- \alpha^2}{\nu'\sqrt{\nu''}}
\end{equation}
and
%
%
%e1.13 #&#
\begin{eqnarray}
\label{edocap1}E_\infty^{\pm}:= \frac{2\nu'\sqrt{\nu''} \pm
\sqrt{4\nu''\nu'^2 - (\nu
''+\nu')(2(\nu''-\nu'+\nu'^2)-\alpha^2\log{{\nu''}/{\nu
'}})}}{\nu
'+\nu''}.\hspace*{-40pt}
\end{eqnarray}

Note that
%
%
%e1.14 #&#
\begin{equation}
\label{seqequalities} E_\infty^{-} \leq E'_\infty
\leq E_\infty.
\end{equation}
Furthermore, $E'_\infty= E_\infty$ if and only if $E_\infty=
E_\infty
^{-}$ if and only if $\alpha^2 = 0$; that is, any equality in \eqref
{seqequalities} implies a triple equality. It occurs if and only if the
mixture is pure; see \eqref{eq2p}.

%
%pr1 #&#
\begin{proposition}\label{remarkmaxima} For any mixture $\nu$ and any
$k\geq0$, the $k$-complexity functions $\theta_{k, \nu}(u)$ satisfy the
following:

\begin{longlist}[(1)]
\item[(1)] $\theta_{k, \nu}(u)$ is continuous on $\R$ and
differentiable on $\R\setminus\{-E_\infty\}$.

\item[(2)] $\theta_{k, \nu}(u)$ is strictly increasing on $(-\infty
,-E_\infty')$ and strictly decreasing on $(-E_\infty',\infty)$. Its
unique maximum is independent of $k$ and equal to
%
%
%e1.15 #&#
\begin{equation}
\label{totalmaxima} \Sigma_\nu:= \theta_{k,\nu}
\bigl(-E_{\infty}'\bigr) = \frac{1}{2} \log
\frac{\nu
''}{\nu'}-\frac{\nu''-\nu'}{\nu''+ \nu'} > 0.
\end{equation}
\item[(3)] $\theta_{k, \nu}(u)$ has exactly two distinct zeros. The
largest zero is given by $-E_\infty^{-}$ and therefore is independent
of $k$.
\item[(4)] For any $k,k' \geq0$ with $k < k'$, $\theta_{k,\nu}(u) >
\theta_{k',\nu}(u) $ for all $u\in(-\infty,-E_\infty)$.
\item[(5)] For any $k,k' \geq0$ with $k < k'$, $\theta_{k,\nu}(u) =
\theta_{k',\nu}(u) $ for all $u\in[-E_\infty, \infty)$.
\end{longlist}
\end{proposition}

From Theorem \ref{maintheorem} and Proposition \ref{remarkmaxima} we obtain:
%
%
%co1.1 #&#
\begin{corollary}\label{coracaocor}
The mean total number of critical points of index $k$ satisfies
%
%
%e1.16 #&#
\begin{equation}
\lim_{N\to\infty} \frac{1}{N} \log\E\Crt_{N,k}(\R) =
\Sigma_\nu.
\end{equation}
Furthermore, if $B = (-\infty,u)$ with $u \leq-E'_\infty$, then
%
%
%e1.17 #&#
\begin{equation}
\lim_{N\to\infty} \frac{1}{N} \log\E\Crt_{N,k}(-
\infty,u) = \theta_{k,\nu}(u).
\end{equation}
\end{corollary}

%
%re1 #&#
\begin{remark}
By symmetry, Theorem \ref{critical2} also holds as stated for the
random variables $\Crt_{N,N-l}(B)$, with $l \geq1$ fixed if one
replaces $\theta_{k,\nu}(u)$ by $\theta_{k,\nu}(-u)$.
\end{remark}

We now use Theorem \ref{maintheorem} and Proposition \ref{remarkmaxima}
to describe the bottom landscape of~$H_{N}$. For any integer $k \geq
0$, we introduce $E_k = E_k(\nu)>0$ as the unique solution in
$(E_\infty
,\infty)$ to (see Figure~\ref{fig1} again)
%
%
%e1.18 #&#
\begin{equation}
\label{defEk} \theta_{k,\nu}\bigl(-E_k(\nu)\bigr) = 0.
\end{equation}
That is, $-E_k(\nu)$ is the smallest zero of the $k$-complexity function.
It is important to note that, by items (4) and (5) of Proposition \ref
{remarkmaxima}, the sequence $(E_k(\nu))_{k\in\N}$ is nonincreasing.
Its structure is of extreme importance and will be also explored in
Section~\ref{secSP}. We have the following consequence of Theorem~\ref{maintheorem}:

%
%th1.2 #&#
\begin{theorem} \label{tnofiniteindex}
For $k\geq0$ and $\varepsilon> 0$, let $A_{N,k}(\varepsilon)$ be the
event ``there is a critical value of $H_{N}$ below the
level $-N(E_{k}(\nu)+\varepsilon)$ and with index larger or equal to $k$,''
that is,
\[
A_{N,k}(\varepsilon)=\Biggl\{\sum_{i=k}^\infty
\Crt_{N,i}\bigl(\bigl(-\infty,-E_k(\nu)-\varepsilon\bigr)
\bigr)>0\Biggr\}
\]
and
$B_{N,k}(\varepsilon)$ be the
event ``there is a critical value of index $k$ of $H_{N}$
above the level $-N(E_\infty^{-}-\varepsilon)$,'' that is,
\[
B_{N,k}(\varepsilon)=\bigl\{\Crt_{N,k}\bigl(
\bigl(-E_\infty^{-}+\varepsilon,\infty\bigr)\bigr)>0\bigr\}.
\]
Then for all $k\ge0$ and $\varepsilon>0$,
%
%
%e1.19 #&#
\begin{equation}\qquad
\limsup_{N\rightarrow\infty} \frac{1}{N} \log\Pro\bigl(A_{N,k}(
\varepsilon)\bigr) < 0 \quad\mbox{and} \quad\limsup_{N\rightarrow\infty}
\frac{1}{N} \log\Pro\bigl(B_{N,k}(\varepsilon)\bigr) < 0.
\end{equation}
\end{theorem}

Theorem \ref{tnofiniteindex} says that with overwhelming probability
all critical values of $H_{N}$ of index $k$ are inside the interval
$[-NE_{k}, -NE_\infty^{-}]$. A similar result was derived for the pure
case in \cite{ABC}.
However, in the pure case it was shown (Theorem 2.2 of~\cite{ABC}) that
the probability of finding a critical point of finite index above the
level $-N E_\infty$ is asymptotically of order $\exp(-N^2 C)$.

We now study the number of critical points with diverging index and the
total number of critical points (regardless of index). Let $k=k(N)$ be
a sequence of integers such that as $N$ goes to infinity,
%
%
%e1.20 #&#
\begin{equation}
\label{ek} \frac{k(N)}{N} \rightarrow\gamma\in(0,1).
\end{equation}
Let $s_\gamma\in(-\sqrt{2},\sqrt{2})$ be defined as solution of
%
%
%e1.21 #&#
\begin{equation}
\label{ihgai} \frac{1}{\pi}\int_{-\sqrt{2}}^{-s_{\gamma}}
\sqrt{2-x^2} \,\d x = \gamma.
\end{equation}

Our next result is the analogue of Theorem \ref{maintheorem} for
critical points of diverging index.
%
%
%th1.3 #&#
\begin{theorem}\label{midguythm} For any sequence $k(N)$ satisfying
\eqref{ek}, as $N$ goes to infinity,
\begin{eqnarray*}
% \begin{eqnarray}
&&\lim_{N\rightarrow\infty} \frac{1}{N} \log\E
\Crt_{N,k(N)}(B)
\\
&&\qquad= \sup_{y \in B} \biggl\{ \frac{1}{2} \log
\frac{\nu''}{\nu'} + \frac
{1}{2} \biggl(s_\gamma^2 -
\frac{2\nu''}{\alpha^2} \biggl(s_\gamma-\frac{\nu
'y}{(2\nu'')^{{1}/{2}}}
\biggr)^2 - y^2 \biggr) \biggr\}
\\
&&\qquad:= \sup_{y \in B} \theta_{\gamma,\nu}(u).
\end{eqnarray*}
\end{theorem}

%
%re2 #&#
\begin{remark}
From Theorem \ref{midguythm} one can easily get analogues of Theorem
\ref{tnofiniteindex} and Corollary \ref{coracaocor} for the case of
critical points with diverging index. Its statements are adapted
rewrites of the respective results. We leave this to the reader.
\end{remark}

We also provide the complexity for the expected total number of
critical values at a level of energy. Precisely, define
%
%
%e1.22 #&#
\begin{equation}
\label{eqqsela} \theta_\nu(u) = \cases{ \theta_{0,\nu}(u) &\quad $
\mbox{if } u\leq-E_\infty',$ \vspace*{2pt}
\cr
\theta_{0,\nu}(-u) &\quad$\mbox{if } u\geq E_\infty',$
\vspace*{2pt}
\cr
\displaystyle\frac{1}{2} \biggl(\log\frac{\nu''}{\nu'} -
\frac{\nu''-\nu'}{\nu
'^2-\nu'+\nu''}u^2 \biggr) &\vspace*{2pt}\cr
\qquad= \sup_{\gamma\in(0,1)}
\theta_{\gamma,\nu}(u), &\quad $\mbox{otherwise}.$}
\end{equation}
%
%
%th1.4 #&#
\begin{theorem}
The total number of critical points satisfies
\label{tcomplexityglobal}
%
%
%e1.23 #&#
\begin{equation}
\label{totalequation} \lim_{N\rightarrow\infty} \frac{1}{N} \log
\E
\Crt_{N}(B) = \sup_{u\in B} \theta_\nu(u)
:= \Theta_\nu(u).
\end{equation}
\end{theorem}

%
%re3 #&#
\begin{remark}
The last result can be interpreted as follows: the mean number of
critical points at levels of the form $Nu +o(N)$ is asymptotically
given by the mean number of local minima, local maxima or critical
points of index $k(N) \sim\gamma(u) N$ if $u\leq-E_{\infty}', u\geq
E_{\infty}', -E_{\infty}'\leq u \leq E_{\infty}'$, respectively. Here,
$\gamma(u) \in(0,1)$ is such that $s_{\gamma(u)}= \sqrt{2} \frac
{u}{E_{\infty}'}$; see \eqref{ihgai}.
\end{remark}

We also investigate the landscape of the Hamiltonian $H_{N}$ by
analyzing the mean Euler characteristic of level sets as $N$ goes to
infinity. In order to state our results we need further notation. The
Hermite functions $\phi_j$, $j\in\mathbb N$, are defined by
%
%
%e1.24 #&#
\begin{equation}
\label{HFd} \phi_j(x) = \bigl(2^jj!\sqrt{\pi}
\bigr)^{-1/2} h_j(x) e^{-{x^2}/{2}},
\end{equation}
where $h_j$, $j\in\mathbb N$ are Hermite polynomials,
%
%
%e1.25 #&#
\begin{equation}
\label{eq1} h_j(x) = e^{x^2}\biggl(-\frac{\d}{\d x}
\biggr)^j e^{-x^2}.
\end{equation}
In particular, $h_0(x) =1, h_1(x)=2x, h_2(x)=4x^2-2x. $ The Hermite
functions are orthonormal functions in $\R$ with respect to Lebesgue measure.

We denote by $\chi(A_u)$ the Euler characteristic of a level set
\[
A_u:= \bigl\{\bolds\sigma\in S^{N-1}(\sqrt{N})\dvtx
H_{N}(\bolds\sigma) \leq N u \bigr\}.
\]
$\chi(\cdot)$ is a topological invariant, integer valued function that
is defined for any CW-complex as the alternating sum of Betti's numbers
\cite{Betti}. It is a functional that is invariant under homotopies and
satisfies
%
%
%e1.26 #&#
\begin{eqnarray}
\label{eqeulde} \chi(A \cup B)& =& \chi(A)+\chi(B)-\chi(A \cap B),\qquad
\chi(\mathbb
B)=1 \quad\mbox{and}
\nonumber
\\[-8pt]
\\[-8pt]
\nonumber
 \chi(S_N) &=& 1 + (-1)^{N-1},
\end{eqnarray}
where $\mathbb B$ denotes a $N$-dimensional unit ball, $S_N$ the
$N$-dimensional unit sphere and $A$, $B$ are CW-complexes. $\chi(\cdot
)$ roughly measures the number of connected components and its number
of attached cylindrical holes and handles. Since we are only interested
in Euler characteristics of level sets of functions that are almost
surely Morse, we use the equivalent definition that follows from
Morse's theorem (see \cite{AT07}, Theorem 9.3.2),
\[
\chi(A_u):= \sum_{k=0}^{N-1}
(-1)^{k} \Crt_k(A_u).
\]

The strategy of using Rice's formula to compute Euler characteristics
of level sets was developed in \cite{AT07,TaylorJon,Taylor22} and also
explored in \cite{Azaisbook}. In fact, in a similar fashion, we prove
the following proposition:

%
%pr2 #&#
\begin{proposition}\label{eulerexact}
%
%
%e1.27 #&#
\begin{eqnarray}\label{eqmix}
&&\E\chi(A_u)
\nonumber\\
&&\qquad=(-1)^{N-1} \biggl(\frac{\nu''}{\nu'} \biggr)^{{(N-1)}/{2}}
\frac{2^{-(N-1)}N}{\sqrt{\pi}\Gamma({N}/{2})}\\
&&\quad\qquad{}\times\int_{-\infty
}^{\infty} \int
_{-\infty}^{u} h_{N-1} \biggl(
\frac{\sqrt{N}(\nu'x -\alpha
y)}{\sqrt{2\nu''}} \biggr) e^{-{N}/{2} (x^2+y^2)} \,\d x \,\d
y.\nonumber
\end{eqnarray}
\end{proposition}

This allows us to derive the asymptotic formula for $\E\chi(A_u)$ and
its relation to the asymptotic complexity of the total number of
critical points; see \eqref{totalequation}.

%
%th1.5 #&#
\begin{theorem}\label{meaneulerasym} The mean Euler--Poincar\'{e}
characteristic $\E\chi(A_u)$ satisfies the following:

\begin{longlist}[(3)]
\item[(1)] If $u \leq-E_\infty'$,
%
%
%e1.28 #&#
\begin{equation}
\E\chi(A_u) = C(N,\nu,u)N^{-{1}/{2}} e^{N \Theta_\nu(u)}\bigl
(1 +
O\bigl(N^{-1}\bigr)\bigr),
\end{equation}
where $C(N, \nu, u)$ is a positive constant given in \eqref{eqCt42}.

\item[(2)] If $ -E_\infty' < u \leq0 $, with $u = -E_\infty' \cos
\omega$, $\omega\in(0,\pi)$
%
%
%e1.29 #&#
\begin{eqnarray}
\E\chi(A_u) &=& (-1)^{N-1}\frac{c(N,\nu)}{2^{{1}/{4}}\pi^{
{1}/{2}} N^{{5}/{4}}}
\frac{e^{N \Theta_\nu(u)}}{f(\omega)(\sin
\omega)^{1/2}} \sin\bigl[ N \tau(\omega) + \rho(\omega) \bigr]
\nonumber
\\[-8pt]
\\[-8pt]
\nonumber
&&{}\times \bigl(1 + O\bigl(N^{-1}\bigr)\bigr),
\end{eqnarray}
where
\[
\tau(\omega) = \frac{1}{2} (\sin2\omega- 2\omega),\qquad \rho(\omega)
= -
\frac{1}{2}\tau(\omega) + \frac{3\pi}{4} + \alpha(\omega),
\]
$c(N,\nu)$ is given in \eqref{cnnu} and $f(\omega)$, $\alpha(\omega)$
are given in \eqref{finalalpha}.
\item[(3)] If $u>0$, we have $\E\chi(A_u) = \E\chi(A_{-u})$ for $N$
even and $\E\chi(A_u) = 2-\E\chi(A_{-u})$ for $N$ odd.
\end{longlist}
\end{theorem}

Let us now describe in words the landscape picture emerging from
Theorem \ref{meaneulerasym}.
Roughly speaking, Theorem \ref{meaneulerasym} says that the mean Euler
characteristic of $A_u$ is in absolute value asymptotically equal to
the total number of critical points \textit{at level} $Nu$ if $u <
E_0$. This picture is fairly intuitive and easy to explain in the
bottom of the landscape. As we increase the energy level $u$ from
negative infinity to $- E_\infty'$, the level set $A_u$ is
``essentially'' a union of disjoint simply connected neighborhoods of
local minima. Since these are exponentially large and dominate the
total number of critical points, the mean Euler characteristic is
positive and of the same size.
As we cross the level $-E_\infty'$, local minima cease to dominate. The
total number of critical points and the Euler characteristic (in
absolute value) is given by the critical values of dominant divergent
index. The landscape is then hard to visualize. By increasing a tiny
amount of energy it oscillates from a large positive to a large
negative Euler characteristic (and vice versa). This oscillation
continues up to level $ E_\infty'$. It would be of interest to find a
simple and intuitive geometric reason for this large oscillation. By
symmetry above $ E_\infty'$ we have ``essentially'' covered the whole
sphere minus an exponentially large number of disjoint simply connected sets.

The rest of the paper is organized as follows. In Section~\ref{sec2} we
prove all Theorems about the complexity function. Their proofs follow
the same strategy of \cite{ABC}. Namely, they will follow from an exact
formula for the mean number of critical points of index $k$ that
translates the problem to a Random Matrix Theory question. This formula
is more involved than the pure case since in a mixture the Hessian
matrix gains an independent Gaussian component on the diagonal. This
leads to a different variational principle that we analyze. In
Section~\ref{eulerproof} we prove the results related to the Euler's
characteristic. In Sections~\ref{secSP} and \ref{sec5} we explain our interest in
such functions, and we relate $H_N$ to Hamiltionians of classical
models in statistical physics.

%s2 #&#
\section{Complexity of critical points}\label{sec2}

%s2.1 #&#
\subsection{Main identity}
In this section, we introduce the main identity that relates the mean
number of critical points of index $k$ with the $k$th smallest
eigenvalue of the Gaussian orthogonal ensemble. This identity, given in
Proposition \ref{identity}, is the analogous of Theorem 2.1 of \cite
{ABC} and it is the first step of the proofs of Theorems~\ref
{critical2}, \ref{tnofiniteindex}, \ref{tcomplexityglobal} and Proposition~\ref{theoasd}.

We fix our notation for the Gaussian orthogonal ensemble (GOE). The GOE
is a probability measure on the space of real symmetric
matrices. Namely, it is the probability distribution of the $N \times
N$ real
symmetric random matrix $M^N$, whose entries $(M_{ij}, i\leq j)$ are
independent centered Gaussian random variables with variance
%
%
%e2.1 #&#
\begin{equation}
\label{eMs} \mathbb E M_{ij}^2 = \frac{1+\delta_{ij}}{2N}.
\end{equation}
We will denote by $\E^N_{\GOE}$ the expectation under the GOE
ensemble of
size $N\times N$.

Let $\lambda^N_0 \leq\lambda^N_1\le\cdots\leq\lambda^N_{N-1}$ be
the ordered
eigenvalues of $M^N$.

%
%pr3 #&#
\begin{proposition}\label{identity}
The following identity holds for all $N$, $\nu$,
$k\in\{0,\ldots,\break N-1\}$, and for all open sets $B \subset\mathbb R$:
%
%
%e2.2 #&#
\begin{eqnarray}
\label{eexactk}  &&\E\bigl[\Crt_{N,k}(B)\bigr]\nonumber\hspace*{-35pt}\\
&&\qquad = C\bigl(N,
\nu',\nu''\bigr)\hspace*{-35pt}
\\
&&\qquad\quad{}\times\int
_{B} \E_\GOE^N \biggl[ \exp\biggl\{
\frac{N}{2} \biggl(\bigl(\lambda_{k}^{N}
\bigr)^2 - y^2 - \frac{2\nu''}{\alpha^2} \biggl(
\lambda_{k}^{N}-\frac{\nu'y}{(2\nu
'')^{{1}/{2}}} \biggr)^2
\biggr) \biggr\} \biggr]\,\d y,\nonumber\hspace*{-35pt}
\end{eqnarray}
where $C(N,\nu',\nu'')= 2 \sqrt{\frac{2\nu''N}{\nu'\pi\alpha
^2}}(\frac
{\nu''}{\nu'})^{ N/2} \frac{\nu'}{\sqrt{2\nu''}}$.
\end{proposition}

\begin{pf}
Proof of Proposition \ref{identity} is a rewrite of the proof of
Theorem 2.1 of~\cite{ABC} with one subtle difference: the law of the
Hessian in the mixed case gains an independent Gaussian component on
its diagonal. In this proof, we use $H$ to denote~$H_{N}$.

The hypothesis on $\nu$ allows us to apply Rice's formula, in the form
of Lemma~3.1 of \cite{ABC}. It says that
using $\d\bolds\sigma$ to denote the usual surface measure on
$S^{N-1}(\sqrt{N})$,
%
%
%e2.3 #&#
\begin{eqnarray}
\label{emetak}&& \E\Crt_{N,k}(B)\nonumber\hspace*{-35pt}\\
&&\qquad= \int_{S^{N-1}(\sqrt{N})} \E
\bigl[ \bigl| \det\nabla^2H(\bolds\sigma) \bigr| \indi\bigl\{
H(\bolds
\sigma) \in N B,i\bigl(\nabla^2 H(\bolds\sigma)\bigr)=k\bigr\}|\hspace*{-35pt}\\
&&\hspace*{214pt}\quad\qquad \nabla H(\bolds\sigma) = 0 \bigr]
\phi_{\bolds\sigma}(0) \,\d\bolds\sigma,\nonumber\hspace*{-35pt}
\end{eqnarray}
where $\phi_{\bolds\sigma}$ is the density of the gradient vector
of $H$.

Now, since $H$ is invariant under rotations, to compute the above
expectation it is enough to study the joint distribution of $(H,\nabla
H, \nabla^2 H)$ at the north pole $\mathbf{n}$. We fix a
orthogonal base for the tangent plane at the north pole, and we
consider $\nabla H (\mathbf{n}), \nabla^2 H(\mathbf{n})$ with
respect to that base. Denoting subscript by a derivative according to a
orthonormal basis in $T_{\bold\sigma} S^{N-1}(\sqrt{N})$ we have that

%
%le1 #&#
\begin{lemma}
\label{lconditioning}
For
all $1\le i\le j\le N-1$,
\begin{eqnarray*}
\label{ecovariancesa} %
\mathbb E\bigl[H(\mathbf{n})^2
\bigr]&=& N,\qquad\mathbb E\bigl[H(\mathbf{n})H_i(\mathbf{n})\bigr]= \mathbb
E\bigl[H_i(\mathbf{n})H_{jk}(\mathbf{n})\bigr]=0,
\\
\mathbb E\bigl[H(\mathbf{n}) H_{ij}(\mathbf{n})\bigr]&=&-
\nu' \delta_{ij},\qquad\mathbb E\bigl[H_i(\mathbf{n})H_j(\mathbf{n})
\bigr]=\nu' \delta_{ij}
\end{eqnarray*}
and
\[
\label{ecovariancesb} \mathbb E\bigl[H_{ij}(\mathbf{n})H_{kl}(\mathbf{n})\bigr]= \frac{1}{N}\bigl[
\nu''(\delta_{ik}\delta_{jl}+
\delta_{il}\delta_{jk})+ \bigl(\nu''
+ \nu'\bigr) \delta_{ij}\delta_{kl}\bigr].
\]
Furthermore, under the conditional distribution $\mathbb P[\cdot|
H(\mathbf{n})=x]$ the random variables $H_{ij}(\mathbf{n})$ are
Gaussian variables with
\[
\mathbb E\bigl[H_{ij}(\mathbf{n})\bigr]=- \frac{x}{N}
\nu'\delta_{ij}
\]
and
\[
\mathbb E \bigl[H_{ij}(\mathbf{n})H_{kl}(\mathbf{n})
\bigr]= \frac
{1}{N} \bigl[\nu''(1+
\delta_{ij})\delta_{ik}\delta_{jl} +
\alpha^2 \delta_{ij}\delta_{kl}\bigr],
\]
that is, if $M^{N-1}$ is distributed as a $(N-1)\times(N-1)$ GOE matrix
\begin{eqnarray*}
\E\bigl[\nabla^2 H | H(\mathbf{n})\bigr] &\stackrel{d} {=}&\biggl(
\frac{N-1}{N}2\nu''\biggr)^{1/2}M^{N-1}\\
&&{}+ \frac{1}{\sqrt{N}}\biggl(\alpha Z - \frac{1}{\sqrt{N}}\nu' H(
\mathbf{n})\biggr) I,
\end{eqnarray*}
where $Z$ is an independent standard Gaussian.
\end{lemma}

The above lemma implies that \eqref{emetak} can be rewritten as
%
%
%e2.4 #&#
\begin{eqnarray}
\label{ebaa} &&\mathbb E \Crt_{N,k}(B)
\nonumber\\
&&\qquad= \omega_N \E\biggl[ \E\biggl[ \biggl\llvert\det\biggl(
\biggl(
\frac{N-1}{N}2\nu''\biggr)^{1/2}M^{N-1}
+ \frac{1}{N}\bigl(\sqrt{N}\alpha Z - \nu' H(\mathbf{n})
\bigr) I \biggr)\biggr\rrvert
\nonumber\\
&&\hspace*{38pt}\qquad\quad{}\times\indi\biggl\{i \biggl[\biggl(\frac{N-1}{N}2\nu''
\biggr)^{{1}/{2}}M^{N-1} + \biggl(\alpha\frac{Z}{\sqrt{N}} -
\nu'\frac{H(\mathbf{n})}{N}\biggr) I \biggr]=k \biggr\}\\
&&\hspace*{203pt}\qquad\quad{}\times\indi
\bigl\{H(
\mathbf{n}) \in N B\bigr\} | H(\mathbf{n}) \biggr] \biggr]\nonumber\\
&&\qquad\quad{}\times
\phi_{\mathbf{n} }(\mathbf{n}),\nonumber
\end{eqnarray}
where $\omega_N$, the volume of the sphere $S^{N-1}(\sqrt{N})$ and
$\phi
_{\mathbf{n}}(\mathbf{n})$ are given by
%
%
%e2.5 #&#
\begin{equation}
\label{eaab} \omega_N=(\sqrt{N})^{N-1}\frac{2\pi^{N/2}}{\Gamma(N/2)},\qquad
\phi_{\mathbf{n}}(\mathbf{n})= \bigl(2 \pi\nu'
\bigr)^{-(N-1)/2}.
\end{equation}

Since we can assume $\alpha\neq0$ (the case $\alpha=0$, that is, the
pure p-spin was treated in \cite{ABC}), we can rewrite the conditional
expectation in \eqref{ebaa} as
%
%
%e2.6 #&#
\begin{eqnarray}
\label{poiloas}&& \frac{\sqrt{N}}{\sqrt{2\pi}} \biggl(2 \nu''
\frac{N-1}{N}\biggr)^{{(N-1)}/{2}}
\nonumber\hspace*{-35pt}
\\[-8pt]
\\[-8pt]
\nonumber
&&\qquad{}\times \int_{B}
e^{{-N y^2}/{2}} \E\bigl\llvert\det\bigl(M^{N-1} - X(y) \bigr
) I \bigr
\rrvert\indi\bigl\{i \bigl[M^{N-1} - X(y) I \bigr]=k \bigr\} \,\d y,\hspace*{-35pt}
\end{eqnarray}
where $X(y)$ is a Gaussian random variable with mean $m= \frac{\sqrt
{N}\nu'y}{(2\nu''(N-1))^{1/2}}$ and variance $t^2 = \frac{\alpha
^2}{2\nu
''(N-1)}$.
Hence, we can apply Lemma 3.3 of \cite{ABC} with $G=\R$ to get that~\eqref{poiloas} is equal to
%
%
%e2.7 #&#
\begin{eqnarray}
\label{elosi}\qquad  &&\frac{\Gamma({N}/{2})({(N-1)}/{N})^{-
{N}/{2}}}{\sqrt{\pi t^2}}
\nonumber
\\[-8pt]
\\[-8pt]
\nonumber
&&\qquad{}\times \int_{B}
\E_\GOE^N [ \exp\biggl\{\frac{N}{2} \biggl(\bigl(
\lambda_{k}^{N}\bigr)^2 - y^2 -
\frac{2\nu''}{\alpha^2} \biggl(\lambda_{k}^{N}-\frac{\nu
'y}{(2\nu
'')^{{1}/{2}}}
\biggr)^2 \biggr) \biggr\} \,\d y.
\end{eqnarray}
Putting \eqref{ebaa}, \eqref{eaab} and \eqref{elosi} together, we
end the proof of Proposition \ref{identity}.
\end{pf}

%s2.2 #&#
\subsection{\texorpdfstring{Proof of Theorems \protect\ref{critical2}, \protect\ref{tnofiniteindex},
\protect\ref{midguythm} and \protect\ref{tcomplexityglobal}}
{Proof of Theorems 1.1, 1.2, 1.3 and 1.4}}\label{sectionlo}

%s2.2.1 #&#
\subsubsection{\texorpdfstring{Proving Theorem \protect\ref{critical2} and Proposition \protect\ref{remarkmaxima}}
{Proving Theorem 1.1 and Proposition 1}}

In this subsection, we will compute the logarithm asymptotics of the
left-hand side of \eqref{eexactk}.

Let $F\dvtx\R^2 \rightarrow\R$ be given by
%
%
%e2.8 #&#
\begin{equation}
\label{fdefini} F(\lambda,y) = \frac{1}{2} \biggl( - \frac{\nu
''+\nu'}{\nu''+\nu'-\nu
'^2}y^2
+ \frac{2\sqrt{2}\sqrt{\nu''}\nu'}{\nu''+\nu'-\nu'^2} \lambda
y-\frac{\nu''-\nu'+\nu'^2}{\nu''+\nu'-\nu'^2} \lambda^2 \biggr).\hspace*{-35pt}
\end{equation}
Note that $ F(\lambda,y) = -a y^2+ by\lambda-c \lambda^2$ for some
constants $a,b,c>0$. Let
%
%
%e2.9 #&#
\begin{eqnarray}
I_1(x) &=& \int_{\sqrt{2}}^{x}
\sqrt{z^2-2} \,\d z
\nonumber
\\[-8pt]
\\[-8pt]
\nonumber
&=& \frac
{1}{2} \bigl(x \sqrt{x^2-2}+
\log[2]-2 \log\bigl[ \bigl(x+\sqrt{x^2-2} \bigr) \bigr] \bigr).
\end{eqnarray}

For any $k \in\N$ fixed, let
%
%
%e2.10 #&#
\begin{equation}
\label{complexityfunction2} \theta_{k,\nu}(u) = %
\cases{\displaystyle
\frac{1}{2} \log\frac{\nu''}{\nu'} + F(-\sqrt{2},u), \qquad \mbox{if }
-E_\infty\leq u, \vspace*{2pt}
\cr
\displaystyle\frac{1}{2} \log
\frac{\nu''}{\nu'} + F\bigl(\lambda^*_k[u],u\bigr)-(k+1)
I_1\bigl(\bigl|\lambda^*_k[u]\bigr|\bigr),\vspace*{2pt}\cr\hspace*{131pt} \mbox{if } u
\leq-E_\infty,} %
\end{equation}
where $ \frac{\nu'\sqrt{2\nu''}u}{\nu''-\nu'+\nu'^2}<\lambda
^*_k[u]\leq
-\sqrt{2}$ is given by
\[
\Psi'\bigl(\lambda^*_k[u]\bigr) = 0, \qquad \Psi(x)=
\frac{2 \nu' \sqrt{2\nu
''}}{\alpha^2} u x- \frac{\nu''-\nu'+\nu'^2}{\alpha^2} x^2 - 2(k+1)
I_1\bigl(|x|\bigr),
\]
that is, $\lambda^*_k[u]$ is a solution on $(-\infty,-\sqrt{2}]$ of
%
%
%e2.11 #&#
\begin{equation}
\label{inlambda} \frac{\nu' \sqrt{2\nu''}}{\alpha^2} u - \frac
{\nu''-\nu'+\nu'^2}{\alpha
^2} \lambda^*_k[u]
+ (k+1) \sqrt{\bigl(\lambda^*_k[u]\bigr)^2-2}= 0.
\end{equation}

Our goal in this section is to prove that $\theta_{k,\nu}$ is the
$k$-complexity function. When $k=0$ the formula for $\theta_{0,\nu}$
simplifies as follows.
%
%
%pr4 #&#
\begin{proposition}\label{equalityinthe} For all $u \in\R$,
%
%
%e2.12 #&#
\begin{equation}
\qquad\theta_{0,\nu}(u)
= %
\cases{ \displaystyle\frac{1}{2} \biggl(\log
\biggl[\frac{\nu''}{\nu'}\biggr] -\frac{u^2 (\nu'+\nu'')}{\nu
'-\nu'^2+\nu''}+\frac{4 u \nu' \sqrt{\nu''}}{\nu'-\nu'^2+\nu''}\vspace*{2pt}\cr
\hspace*{106pt}{}-
\displaystyle\frac{2
(-\nu'+\nu'^2+\nu'' )}{\nu'-\nu'^2+\nu''} \biggr), \vspace*{2pt}\cr
\qquad\mbox
{if } -E_\infty\leq u,
\vspace*{2pt}
\cr
\displaystyle\frac{1}{2}\log\bigl[\nu'-1\bigr]-
\frac{u^2 (\nu'-2)}{4 (\nu'-1)}- I_1\biggl(-\frac{u
\nu'}{\sqrt{2} \sqrt{\nu'(\nu'-1)}}\biggr), \vspace*{2pt}\cr
\qquad\mbox{if } u \leq-E_\infty.} %
\end{equation}
\end{proposition}
%
%
%re4 #&#
\begin{remark}\label{abccase}
It is possible to recover all complexity functions of the pure case by
taking $\alpha$ to zero (i.e., recover the first results of \cite
{ABC}). In particular, if $\alpha=0$, $E'_\infty= E_\infty$, and we do
not have the intermediate regions where the $k$-complexity functions
are equal for different $k$ and nonconstant.
\end{remark}

We postpone the proof of Proposition \ref{equalityinthe} to the end of
this subsection since we will need another characterization of $\theta
_{k, \nu}$.

\begin{pf*}{Proof of Theorem \ref{critical2}}
To prove Theorem \ref{critical2} it suffices to show that $\theta
_{k,\nu
}(u)$ is the logarithm asymptotic limit of the left-hand side of \eqref
{eexactk}.

First, note that we can rewrite \eqref{eexactk} as
%
%
%e2.13 #&#
\begin{equation}
\label{acopado} C_N \E e^{-N \Lambda(\lambda^N_k,Y_N)} \indi\{Y_N
\in B
\},
\end{equation}
where $Y_N$ is a Gaussian random variable of mean zero and variance $N$
independent of $\lambda^N_k$, $\E$ is the expectation with respect to
GOE and $Y_N$ and
%
%
%e2.14 #&#
\begin{eqnarray}
\label{ploeamoreadedio} \lim_{N\rightarrow\infty} \frac{1}{N} \log
C_N &=& \frac{1}{2} \log\frac{\nu''}{\nu'},
\nonumber
\\[-8pt]
\\[-8pt]
\nonumber
\Lambda(\lambda,y)= F(\lambda,y) + \frac
{y^2}{2} &=& \frac{1}{2}
\biggl(\lambda^2 - \frac{2\nu''}{\alpha^2} \biggl(\lambda-
\frac{\nu'y}{(2\nu'')^{{1}/{2}}} \biggr)^2 \biggr).
\end{eqnarray}

By the independence of $Y_N$ and $\lambda_k^N$ and Theorem A.1 of
\cite
{ABC}, the sequence of random variables $(\lambda^N_k,Y_N)$ satisfies a
large deviation principle of speed $N$ and rate function
\[
I_k(\lambda,x) = %
\cases{\displaystyle\frac{x^2}{2} + (k+1)
I_1\bigl(|\lambda|\bigr), &\quad $\mbox{if } \lambda\leq- \sqrt{2}$, \vspace*{2pt}
\cr
\infty,& \quad $\mbox{otherwise}.$} %
\]
Therefore, in view of \eqref{acopado} and \eqref{ploeamoreadedio}, we
can apply Laplace--Varadhan lemma (see, e.g.,~\cite{AmirOfer}, Theorem
4.3.1 and
Exercise~4.3.11) and get that
\begin{eqnarray}
\label{erkow} &&\lim_{N\rightarrow\infty} \frac{1}{N} \log\E
\Crt_{N,k}(B)
\nonumber\hspace*{-35pt}
\\[-8pt]
\\[-8pt]
\nonumber
&&\qquad=\frac{1}{2} \biggl[ \log\frac{\nu''}{\nu'} +\max_{x \in B,
\lambda
\leq- \sqrt{2}}
\biggl\{\lambda^2 - \frac{1}{\alpha^2}\bigl(\nu'x-
\sqrt{2\nu''}\lambda\bigr)^2 - 2
I_k(\lambda,x) \biggr\} \biggr].\hspace*{-35pt}
\end{eqnarray}
We will now analyze the above variational principle. We start with the
case of $B=(-\infty, u)$. We want to find
%
%
%e2.15 #&#
\begin{equation}
\label{erkoe}\qquad  \max_{x \leq u, \lambda\leq-\sqrt{2}} \biggl\{-x^2 +
\lambda^2 - \frac
{1}{\alpha^2}\bigl(\nu'x- \sqrt{2
\nu''}\lambda\bigr)^2 -2 (k+1)
I_1\bigl(|\lambda|\bigr) \biggr\}.
\end{equation}

\textit{Case $u\ge- E'_\infty$}:
If $u \ge-E'_\infty$, then we maximize \eqref{erkoe} in $x$ first. The
maximum is obtained at $x = x_\lambda:=\frac{\nu' \sqrt{2\nu
''}}{\nu
''+\nu'} \lambda\leq u$. Plugging $x_\lambda$ back in \eqref{erkoe},
we get an increasing function in $\lambda$, since $I_1(|\lambda|)$ is
itself decreasing. Thus the maximum is realized at
\[
x = x_\lambda, \qquad\lambda= -\sqrt{2}.
\]
This together with \eqref{erkow} proves Theorem \ref{critical2} in the
case $B=(-\infty,u)$ with $- E'_\infty\leq u$.

\textit{Case $u\le- E'_\infty$}:
In the case $u\le- E'_\infty$, $x_\lambda\leq u$ if and only if $
\lambda\leq\frac{\sqrt{2}u}{E'_\infty}$. Therefore if $x^*$ maximizes
\eqref{erkoe}, then
%
%
%e2.16 #&#
\begin{equation}
\label{casdf} x^*= x_\lambda\Leftrightarrow\lambda\leq\frac{\sqrt
{2}u}{E'_\infty
}\quad
\mbox{and} \quad x^*= u \Leftrightarrow\frac{\sqrt{2}u}{E'_\infty} \leq
\lambda\leq-\sqrt{2}.
\end{equation}

If we plug in the correspondent values of $x$ in each region, we note
that in the first case our function is again increasing in $\lambda$.
Furthermore, since at $\lambda= \frac{\sqrt{2}u}{E'_\infty}$,
$x_\lambda
=u$, we are led to the following variational principle valid in both
cases of~\eqref{casdf}:
%
%
%e2.17 #&#
\begin{eqnarray}
\label{rafuta} &&\max_{{\sqrt{2}u}/{E'_\infty} \leq\lambda
\leq-\sqrt{2}} \biggl\{ -u^2 +
\lambda^2 - \frac{1}{\alpha^2}\bigl(\nu'u- \sqrt{2
\nu''}\lambda\bigr)^2 - 2(k+1)
I_1\bigl(|\lambda|\bigr) \biggr\}
\nonumber
\\
&&\qquad= -\biggl(1 + \frac{\nu'^2}{\alpha^2}\biggr)u^2 + \max
_{{\sqrt{2}u}/{E'_\infty
} \leq\lambda\leq-\sqrt{2}} \biggl\{ \frac{2 \nu' \sqrt{2\nu
''}}{\alpha^2} u \lambda-
\frac{\nu''-\nu'+\nu'^2}{\alpha^2}\lambda^2
\nonumber
\\[-8pt]
\\[-8pt]
\nonumber
&&\hspace*{223pt}\qquad{}- 2(k+1) I_1\bigl(|\lambda
|\bigr) \biggr
\}
\\
&&\qquad= -\biggl(1 + \frac{\nu
'^2}{\alpha^2}\biggr)u^2 + \max
_{{\sqrt{2}u}/{E'_\infty} \leq\lambda
\leq-\sqrt{2}} \Psi(\lambda) = \max_{{\sqrt{2}u}/{E'_\infty
} \leq
\lambda\leq-\sqrt{2}} \Gamma(
\lambda).\nonumber
\end{eqnarray}

Note that $\Psi(\lambda)$ is a parabola $a\lambda^2 + b \lambda, a<0$
plus an increasing function. The critical point of the parabola is
given by
%
%
%e2.18 #&#
\begin{equation}
\label{lambdacritical} \lambda_c = \frac{\nu'\sqrt{2\nu''}u}{\nu
''-\nu'+\nu'^2} \geq- \sqrt{2}
\quad\Longleftrightarrow\quad  u \geq- E_\infty.
\end{equation}
Therefore if $u \geq-E_\infty$, $\Psi$ is an increasing function in
$\lambda$, so its maximum is attained at $\lambda=-\sqrt{2}$. This
proves the theorem in the region $-E_\infty\leq u \leq- E'_\infty$.

If $u < -E_\infty$, equation \eqref{lambdacritical} and the facts that
$\Psi'(-\sqrt{2})<0$ and $\Psi'(\lambda_c)>0$ imply that the
maximum is
taken in the interior of the interval $[\lambda_c,-\sqrt{2}]$ at
$\lambda_k^*[u]$. This completes the proof of the theorem in the case
$B=(-\infty,u)$.

Now, it is easy to extend it to any open set $B$. Let $u^*$ be the
point that realizes the $\sup_{\{u\in B \}} \theta_{k,\nu}(u)$. From
the continuity and uniqueness of a local maxima of $\theta_{k,\nu}$, it
is clear that either $u^* = -E_\infty'$ or $u^*$ is in the boundary of
$B$. Assume without loss of generality that there exists an increasing
sequence $u_n$ in $B$ approaching $u^*$. Since $B$ is open, there exist
$\varepsilon_n >0$ such that
\begin{eqnarray*}
\E\bigl(\Crt_{N,k}(-\infty,u_n) - \Crt_{N,k}(-
\infty, u_n-\varepsilon_n)\bigr) &=& \E
\Crt_{N,k}(u_n-\varepsilon_n,u_n)
\\
&\leq& \E\Crt_{N,k}(B)
\\
&\leq& \E\Crt_{N,k}\bigl(-\infty,u^*\bigr).
\end{eqnarray*}
But since $\theta_{k,\nu}$ is continuous and increasing for $u \leq
-E_\infty'$, the above equation implies
\[
\theta_{k,\nu} (u_n) \leq\lim_{N\rightarrow\infty}
\frac{1}{N} \log\E\Crt_{N,k}(B) \leq\theta_{k, \nu}
\bigl(u^*\bigr)
\]
for all $n$,
which proves Theorem \ref{critical2} for any $B$ open.
\end{pf*}

It remains to prove Proposition \ref{equalityinthe}. We first need the
following miraculous lemma.

%
%le2 #&#
\begin{lemma}\label{surla} For all $u<-E_\infty$,
\[
\frac{\partial}{\partial\nu''}\theta_{0,\nu}(u) = 0.
\]
\end{lemma}

\begin{pf} The proof relies on how we derived $\theta_{0,\nu}(u)$.
When $u<-E_\infty$, $\theta_{0,\nu}(u)$ is the maximum over $\lambda$
of the functional $\Gamma$ (that depends on $\nu''$) given in \eqref
{rafuta}. Its maximizer $\lambda^*(u)$ is the smallest root of a second
degree polynomial that\vadjust{\goodbreak} can be derived from \eqref{inlambda}. This
second degree equation is given by $A + B \lambda+ C \lambda^2 = 0$ where
%
%
%e2.19 #&#
\begin{eqnarray}
\label{sasawqwq} A&=& 2+\frac{2 u^2 \nu'^2 \nu''}{ (\nu'-\nu
'^2+\nu''
)^2},
\nonumber
\\
B&=& -\frac{2 \sqrt{2} u \nu' \sqrt{\nu''} ((-1+\nu') \nu'+\nu'')}{
(\nu'-\nu'^2+\nu'' )^2},
\\
C&=& \frac{2 ((-1+\nu')^2 \nu'^2+\nu''^2 )}{ (\nu'-\nu
'^2+\nu'' )^2}.
\nonumber
\end{eqnarray}
Now the chain rule and the fact that $\lambda^*(u)$ is a maximum imply
that $\frac{\partial}{\partial\nu''}\theta_{0,\nu}(u) = 0$ if and only
if $\frac{\partial}{\partial\nu''} ( \Gamma(\lambda^*(u))
)=0$, and this holds if and only if $ (\frac{\partial}{\partial\nu
''} \Gamma)(\lambda^*(u))=0$. The last condition can be written as
a second degree equation of the form
%
%
%e2.20 #&#
\begin{eqnarray}
\label{sesesewq1} &&\frac{1}{2} \biggl(-\frac{u^2 (-\nu'-\nu'')}{
(\nu'-\nu'^2+\nu
'' )^2}-
\frac{u^2}{\nu'-\nu'^2+\nu''}-\frac{2 \sqrt{2} u \nu' \sqrt{\nu
''} \lambda}{ (\nu'-\nu'^2+\nu'' )^2}\nonumber\\
&&\hspace*{135pt}\qquad{}+\frac{\sqrt{2} u \nu
' \lambda}{\sqrt{\nu''} (\nu'-\nu'^2+\nu'' )} \biggr)
\\
&&\qquad{}+ \frac{1}{2} \biggl( -\frac{\lambda^2}{\nu'-\nu'^2+\nu
''}+\frac{
(-\nu'+\nu'^2+\nu'' ) \lambda^2}{ (\nu'-\nu'^2+\nu''
)^2} \biggr)+
\frac{1}{2 \nu''}=0.\nonumber
\end{eqnarray}
Comparing the coefficients of \eqref{sasawqwq} with \eqref{sesesewq1}
one sees that their ratios are constantly equal to $\frac{1}{4\nu''}$.
This immediately implies that they share the same roots. So $\lambda
^*(u)$ indeed satisfies $ (\frac{\partial}{\partial\nu''} \Gamma
)(\lambda^*(u)) = 0$, and the lemma is proven.
\end{pf}
\begin{pf*}{Proof of Proposition \ref{equalityinthe}}
From Lemma \ref
{surla} we know that for $u<-E_\infty$, $\theta_{k,\nu}$ does not
depend on $\nu''$. By choosing $\nu'' = \nu'^2 - \nu' + \varepsilon$ and
taking $\varepsilon$ to zero we get the desired result. Indeed, when
$\varepsilon$ goes to zero
\[
\lambda^*(u) \rightarrow\frac{u \nu'}{\sqrt{2} \sqrt{(\nu'-1)
\nu'}}, \qquad F\bigl(\lambda^*(u),u\bigr)
\rightarrow\frac{-u^2 (\nu'-2)}{4 (\nu'-1)}.
\]
\upqed\end{pf*}

%s2.2.2 #&#
\subsubsection{\texorpdfstring{Proof of Theorem \protect\ref{tnofiniteindex}}{Proof of Theorem 1.2}}
We want to prove that there are no critical values of index $k$ of
$H_N$ above
$-N(E_{\infty}^{-}-\varepsilon)$. The function $\theta_{k,\nu}$ is
strictly decreasing on $(-E_{\infty}^{-},\infty)$. Using Theorem~\ref
{critical2}, we have
\[
\mathbb E \bigl[ \Crt_{N,k}\bigl(\bigl(-E_{\infty}^{-}+
\varepsilon,\infty\bigr)\bigr) \bigr]\le\exp\bigl\{N \theta
_{k,\nu}
\bigl(-E_{\infty}^{-}+\varepsilon\bigr) +o(N) \bigr\}.
\]
The
constant $-E_{\infty}^{-}$ is defined by $\theta_{k,\nu}(-E_{\infty
}^{-})=0$ for all $k$. Therefore,
$\theta_{k,p}(-E_{k}+\varepsilon)= c(k,\nu,\varepsilon)<0$. An
application of
Markov's inequality as
\[
\Pro\bigl( B_{N,k}(\varepsilon) \bigr) \leq\E\bigl[
\Crt_{N,k}\bigl(-E_{\infty
}^{-}+\varepsilon,\infty
\bigr) \bigr] \leq e^{-Nc(k,\nu,\varepsilon)}
\]
proves Theorem~\ref{tnofiniteindex} for the event $B_{N,k}(\varepsilon)$.
The proof for the event $A_{N,k}(\varepsilon)$ is analogous.

%follows from the definition of critical and full mixture and
%Proposition \ref{remarkmaxima}. The second part follows from Markov's
%Inequality and Theorem \ref{critical2} just as in the proof of Theorem
% \label{goolkk}
% \frac{1}{N} \log\Pro\big[\sum_{i=k}^{N-1}\Crt_{N,k}((-E_\infty^{+}+
% &\leq
% \frac{1}{N} \log\sum_{i=k}^{N-1} \E\big[\Crt_{N,k}((-E_\infty^{+} +
% &\leq\theta_{k,\nu}(-E_\infty^{+} + \varepsilon/2) <0

%s2.2.3 #&#
\subsubsection{\texorpdfstring{Proof of Theorem \protect\ref{midguythm}}{Proof of Theorem 1.3}}

The proof of Theorem \ref{midguythm} follows the same steps as the
proof of Theorem \ref{critical2}. First by Lemma 3.5 of \cite{ABC}, for
any $\varepsilon>0$, there exists a constant $c=c(\gamma,\varepsilon)>0$
such that
\[
\label{concetrationmid}\Pro\bigl(\bigl|\lambda_{k}^N -
s_{\gamma}\bigr|> \varepsilon\bigr)\leq e^{-cN^2}.
\]

Therefore if we use Proposition \ref{identity}, \eqref{ploeamoreadedio}
and the above statement we have that for any $\varepsilon>0, \delta> 0$
there exists constants $c=c(\varepsilon), d=d(\varepsilon)$ such that for $N$
large enough
\begin{eqnarray*}
\label{ubmg} \E\Crt_{N,k}(B) &\le& C_N \int
_{B} e^{{N}/{2} ( F(\lambda
_k^N,y) )}\indi\bigl\{ \lambda_k^N
\in(s_\gamma-\varepsilon, s_\gamma+\varepsilon)\bigr\} +
e^{dN} e^{-cN^2}
\\
&\le& C_N \int_B e^{{N}/{2} \sup_{\lambda\in(s_\gamma-\varepsilon,
s_\gamma+\varepsilon)} \{ F(\lambda,y)\} } \,\d y +
e^{dN} e^{-cN^2}
\\
&\le& C_N e^{{N}/{2} \sup_{\lambda\in(s_\gamma-\varepsilon,
s_\gamma
+\varepsilon), y \in B} \{ F(\lambda,y)\} }(1+\delta) + e^{dN}
e^{-cN^2}.
\end{eqnarray*}
On the other hand we have the lower bound
\begin{eqnarray*}
\E\Crt_{N,k}(B) &\ge &C_N \int_{B}
e^{{N}/{2} ( F(\lambda
_k^N,y) )}\indi\bigl\{ \lambda_k^N
\in(s_\gamma-\varepsilon, s_\gamma+\varepsilon)\bigr\}
\\
&\ge& C_N \int_B e^{{N}/{2} \inf_{\lambda\in(s_\gamma-\varepsilon,
s_\gamma+\varepsilon)} \{ F(\lambda,y) \} } \,\d y
\\
&\ge& C_N e^{{N}/{2} \inf_{\lambda\in(s_\gamma-\varepsilon,
s_\gamma
+\varepsilon)}\{\sup_{ y \in B} \{ F(\lambda,y) \} \} }(1-\delta).
\end{eqnarray*}
Taking $\frac{1}{N}\log$ on both bounds and taking $\varepsilon$ to zero
afterward, we see that
\[
\frac{1}{N}\log\E\Crt_{N,k}(B) = \sup_{y \in B}
\bigl\{ F(s_\gamma, y)\bigr\}.
\]
%
%s2.2.4 #&#
\subsubsection{\texorpdfstring{Proof of Theorem \protect\ref{tcomplexityglobal}}{Proof of Theorem 1.4}}
We now prove the asymptotic limit of the mean number of critical points at some
level of energy.

Since the total number of critical points is greater than the number of
critical points of index $k(N)$ with $k(N)$ satisfying \eqref{ek} for
$\gamma\in[0,1]$ we clearly have the lower bound
%
%
%e2.21 #&#
\begin{equation}
\label{ras} \sup_{\gamma\in[0,1]} \sup_{u\in B}
\theta_{\gamma,\nu}(u) \leq\lim_{N\rightarrow\infty}\frac
{1}{N}\log\E
\Crt_{N}(B).
\end{equation}

For $u\leq-E_\infty'$, taking $\gamma= 0$ (i.e., considering the
complexity of local minima) we get the right-hand side of \eqref
{totalequation}. For $ u \in(-E_\infty', E_\infty')$ the supremum on
$\gamma$ of $\theta_{\gamma,\nu}(u)$ is attained at $\gamma\in(0,1)$
such that $s_\gamma= \frac{\sqrt{2}u}{E_\infty'}$. Plugging this value
back on the left-hand side of \eqref{ras}, we get the right-hand side
of \eqref{totalequation}. Last, for $u\geq E_\infty$, one just needs to
take the complexity of local maxima. This is enough to prove a lower bound.

To show a matching upper bound, we proceed as follows. A sum over $k$
in Proposition \ref{identity} gives us that
\begin{eqnarray*}
&&\E\bigl[\Crt_{N}(B)\bigr] \\
&&\qquad= 2N\sqrt\frac{2}{\nu'} \biggl(
\frac{\nu''}{\nu'}\biggr)^{ N/2} \int_{B} \E
_\GOE^N \int\exp\bigl\{N F(z,y) \bigr\} \,\d y
L_N (\d z),
\end{eqnarray*}
and $L_N$ is the empirical spectral measure of the GOE matrix. The
constant in front the integral gives a constant term $C_\nu$ after the
$\frac{1}{N}\log$ limit.
Furthermore,
%
%
%e2.22 #&#
%e2.23 #&#
\begin{eqnarray}
\label{dskop}&& \int_{B} \E_\GOE^N
\int\exp\bigl\{N F(z,y) \bigr\} \,\d y L_N (\d z)\nonumber\\
&&\qquad \leq N \int
_{B} \sup_{z \in\R} \exp\bigl\{NF(z,y) \bigr
\} \,\d y
\\
&&\qquad \leq N \int_{B} e^{-({N}/{2}) {(\nu'' - \nu')}/{(\nu'^2-\nu
' + \nu
'')}y^2} \,\d y.\nonumber
\end{eqnarray}
So if $B \cap(-E_\infty',E_\infty') \neq\varnothing$, this matches the
right-hand side of $\eqref{totalequation}$. If $B \subseteq(-\infty
,-E_\infty')$, then we can estimate \eqref{dskop} with
\[
N \int_{B} \E_\GOE^N \int\exp\bigl
\{N F(\lambda_0,y) \bigr\}.
\]
Applying $\log$, dividing by $N$ and taking limits we get Theorem \ref
{tcomplexityglobal} from Theorem~\ref{critical2}.

%s3 #&#
\section{\texorpdfstring{Proof of Proposition \protect\ref{eulerexact} and Theorem \protect\ref{meaneulerasym}}
{Proof of Proposition 2 and Theorem 1.5}}\label{eulerproof}
In this section we prove Proposition \ref{eulerexact} and Theorem \ref
{meaneulerasym}.

\begin{pf*}{Proof of Proposition \ref{eulerexact}}
We start with the following identity:
\begin{eqnarray*}
\E \chi(A_u) &=& \sum_{k=0}^{N-1}
(-1)^{k} \Crt_k\bigl(A(u)\bigr)
\\[-1pt]
&=& \sum_{k=0}^{N-1} (-1)^{k}
\int_{S^{N-1}(\sqrt{N})} \E\bigl( \bigl|\det\nabla^2
H_N (\sigma)\bigr| \indi_{\{i(\nabla^2 H_N (\sigma)) =k\}} \indi_{\{
H_N(\sigma) \leq Nu\}} | \\[-1pt]
&&\hspace*{238pt}{}\nabla
H_N(\sigma) = 0 \bigr)
\\[-1pt]
&&\hspace*{91pt}{}\times\phi_{\nabla H_N}(0) \,\d\sigma\\
&=& \bigl(2\nu'\pi
\bigr)^{-{(N-1)}/{2}} \bigl|S^{N-1}(\sqrt{N})\bigr| \frac{1}{\sqrt{2\pi
N}}
\\[-1pt]
&&{}\times \sum_{k=0}^{N-1} \int
_{-\infty}^{Nu} \E\bigl( (-1)^{k} \bigl|\det
\nabla^2 H_N (\sigma)\bigr| \indi_{\{ i(\nabla^2 H_N (\sigma)) =k \}} |
H_N(\sigma) = x \bigr)\\[-1pt]
&&\hspace*{52pt}{}\times e^{-({1}/{2N}) x^2} \,\d x
\\[-1pt]
&=& \bigl(2\nu'\pi\bigr)^{-{(N-1)}/{2}} \frac{2 \pi^{
{N}/{2}}}{\Gamma(
{N}/{2})}N^{{(N-1)}/{2}}\\[-1pt]
&&{}\times
\frac{\sqrt{N}}{\sqrt{2\pi}} \int_{-\infty
}^{u} \E\bigl( \det
\nabla^2 H_N (\sigma) | H_N(\sigma) = Nx
\bigr) e^{-({N}/{2}) x^2} \,\d x
\\[-1pt]
&=& \nu'^{-{(N-1)}/{2}} 2^{-{(N-2)}/{2}} \\[-1pt]
&&{}\times\frac{N^{{N}/{2}}}{\Gamma({N}/{2})} \int
_{-\infty}^{u} \E\bigl( \det\nabla^2
H_N (\sigma) | H_N(\sigma) = Nx \bigr) e^{-({N}/{2}) x^2}
\,\d x.
\end{eqnarray*}

%
%le3 #&#
\begin{lemma}\label{charac} If $M^N$ is a $N\times N$ GOE with variance
$\E M_{ij}^2 = \frac{1+\delta_{ij}}{2N}$, then for any $x\in\R$
\[
\E\det\bigl(M^N-xI\bigr) = 2^{-N} N^{-{N}/{2}}(-1)^N
h_N(\sqrt{N}x),
\]
where $h_N(x)$ is given in \eqref{eq1}.
\end{lemma}
\begin{pf}
The proof, a straight-forward linear algebra exercise, can be found as
Corollary 11.6.3 in \cite{AT07}.
\end{pf}

Now by Lemma \ref{lconditioning},
%
%
%e3.1 #&#
\begin{eqnarray}
\label{Laaa2} \E\chi(A_u)& =& \nu'^{-{(N-1)}/{2}}
2^{-{(N-2)}/{2}} \frac{N^{
{N}/{2}}}{\Gamma({N}/{2})}\frac{\sqrt{N}}{\sqrt{2\pi}}
\nonumber
\\[-2pt]
&&{}\times\int_{-\infty}^{\infty} \int_{-\infty}^{u}
\E\biggl( \det\biggl[ \biggl(\frac
{N-1}{N}2\nu''
\biggr)^{1/2}M^{N-1} + \bigl(\alpha y - \nu'x
\bigr)I \biggr] \biggr)\\[-2pt]
&&\hspace*{33pt}\qquad{}\times e^{-({N}/{2}) x^2} e^{-({N}/{2})y^2}
\,\d x \,\d y.\nonumber
\end{eqnarray}
The double integral becomes
\begin{eqnarray*}
&&\biggl(\frac{N-1}{N}2\nu''\biggr)^{{(N-1)}/{2}}\\[-2pt]
&&\qquad{}\times\int_{-\infty}^{\infty} \int_{-\infty}^{u}
\E\biggl( \det\biggl[M^{N-1} +\biggl(\frac{N-1}{N}2\nu
''\biggr)^{-{1}/{2}}\bigl(\alpha y -
\nu'x\bigr)I \biggr] \biggr) \\[-2pt]
&&\qquad\hspace*{56pt}{}\times e^{-{N}/{2}
(x^2+y^2)} \,\d x \,\d y,
\end{eqnarray*}
which by Lemma \ref{charac} can be rewritten as
%
%
%e3.2 #&#
\begin{eqnarray}
\label{impa1}&& (-1)^{N-1}\biggl(\frac{\nu''}{2N}\biggr)^{{(N-1)}/{2}}
\int_{-\infty}^{\infty} \int_{-\infty}^{u}
h_{N-1} \biggl( \frac{\sqrt{N}(\nu'x -\alpha y)}{\sqrt{2\nu
''}} \biggr)
\nonumber
\\[-8pt]
\\[-8pt]
\nonumber
&&\hspace*{126pt}\qquad{}\times e^{-{N}/{2} (x^2+y^2)} \,\d x \,\d y.
\end{eqnarray}
Combining \eqref{Laaa2} and \eqref{impa1} we get Proposition \ref
{eulerexact}.
\end{pf*}\eject

We will need the following lemma to prove Theorem \ref{meaneulerasym}:

%
%le4 #&#
\begin{lemma}\label{lemmanecessa} Let $a$, $b$ be constants such that
$a>1/2$ and $b\geq0$.
Set
\[
I_N(M) = \int_{M}^{\infty}
\phi_{N-1}(\sqrt{N}x)e^{ax^2+bx} \,\d x.
\]
As $N$ goes to infinity:
\begin{longlist}[(1)]
\item[(1)] If $\sqrt{2}\leq M$, then
$I_N(M) = O(e^{-N(aM^2+bM+I_1(M))})$.
\item[(2)] If $-\sqrt{2}<M < \sqrt{2}$ and if we set $M = \sqrt
{2}\cos
\omega$ with $\varepsilon< \omega< \pi-\varepsilon$, then $I_N(M)$ is
equal to
%
%
%e3.3 #&#
\begin{eqnarray}\qquad
&&\frac{2^{-3/4}\pi^{-1/2}e^{-N(aM^2+bM)}}{N^{{5}/{4}}|m'(2\iota
(M))|(\sin\omega)^{1/2}} \sin\biggl[ \biggl(\frac{N}{2}-\frac{1}{4}
\biggr) (\sin2\omega- 2\omega) + \frac{3\pi}{4} + \alpha(M)
\biggr]
\nonumber
\\[-8pt]
\\[-8pt]
\nonumber
&&\qquad{}\times
\bigl(1 + O\bigl(N^{-1}\bigr)\bigr).
\end{eqnarray}
\item[(3)] If $M\leq-\sqrt{2}$, then $ I_N(M)= LN^{-1/2}e^{-
N\lambda
(a,b,M)}$ where $\lambda(a,b,M)$ is the minimum of $ax^2 + bx +
I_1(-x)$ in $[M,-\sqrt{2}]$ and $L$ is a positive constant that depends
on $a,b$ and $M$ as in \eqref{eqclem82}.
\end{longlist}
\end{lemma}

A few comments before the proof of the above lemma. First, under the
assumption that $a>1/2$ and $b>0$ the major contribution to the
integral in part (2) comes from a small neighborhood of M, instead of
the minimum of $ax^2+bx$. This is due to rapid oscillations of $\phi
_{N-1}$ inside the ``bulk'' $(-\sqrt{2},\sqrt{2})$. Second, in part (3),
the condition that the minimizer of $ax^2 + bx + I_1(-x)$ lies inside
$[M,-\sqrt{2}]$ is similar to the condition on \eqref{inlambda}. This
will lead to the asymptotic Euler's characteristic in the region
$u<-E_\infty'$.

The main tool to prove Lemma \ref{lemmanecessa} is the following
well-known formula for the asymptotics of the Hermite functions, first
proved by Plancherel--Rotach \cite{Prota}. Let
\[
h(x) = \biggl|\frac{x - \sqrt{2}}{x + \sqrt{2}} \biggr|^{1/4} + \biggl|\frac{x +
\sqrt{2}}{x - \sqrt{2}}
\biggr|^{1/4}.
\]

%
%le5 #&#
\begin{lemma}[(Plancherel--Rotach \cite{Prota})] \label{lemPR}There
exists $\delta_0 >0$ such that for any $0<\delta<\delta_0$, the
following asymptotics hold uniformly in each region:
\begin{longlist}[(1)]
\item[(1)]
If $x<-\sqrt{2}-\delta$,
\[
\phi_{N-1}(\sqrt{N}x) = (-1)^{N-1} \frac{e^{-N I_1(-x)}}{\sqrt{4\pi
\sqrt{2N}}}h(x)
\bigl(1 + O\bigl(N^{-1}\bigr)\bigr).
\]
\item[(2)]
If $-\sqrt{2}-\delta<x<-\sqrt{2}+\delta$,
%
%
%e3.4 #&#
%e3.5 #&#
\begin{eqnarray*}
\label{PRA4} &&\phi_{N-1}(\sqrt{N}x)\\
&&\qquad= \frac{(-1)^{N-1}}{(2N)^{1/4}}
\biggl\{ \biggl|
\frac{x-\sqrt{2}}{x+\sqrt{2}} \biggr|^{1/4} \biggl|\frac{3N}{2}I_1(-x)\biggr|^{{1}/{6}}
\Ai\biggl[\biggl(\frac
{3N}{2}I_1(-x)\biggr)^{{2}/{3}}
\varepsilon(x) \biggr] \\
&&\hspace*{60pt}\qquad{}\times\bigl(1+O\bigl(N^{-1}\bigr)\bigr)
\\
&&\hspace*{60pt}\qquad{}- \biggl|\frac{x+\sqrt{2}}{x-\sqrt{2}}\biggr |^{1/4}\biggl |\frac
{3N}{2}I_1(-x)\biggr|^{-{1}/{6}}
\\
&&\hspace*{108pt}\qquad{}\times \Ai' \biggl[\biggl(\frac
{3N}{2}I_1(-x)
\biggr)^{{2}/{3}}\varepsilon(x) \biggr]\bigl(1+O\bigl
(N^{-1}\bigr)
\bigr) \biggr\},
\end{eqnarray*}
where $\Ai(x)$ is the Airy function of first kind,
$\Ai(x)= \frac{2}{\pi} \int_{-\infty}^{\infty}
\cos(\frac{t^3}{3}+ tx ) \,\d t$, and $\varepsilon(x)=\frac
{-x-\sqrt{2}}{|-x-\sqrt{2}|}, x \neq-\sqrt{2}$, $\varepsilon
(-\sqrt{2})=0$ and $Ai'(x)$ is the derivative of $Ai(x)$.

\item[(3)]
If $-\sqrt{2}+\delta< x < \sqrt{2}-\delta$ and if we set $x = \sqrt
{2}\cos\omega$ with $\varepsilon< \omega< \pi-\varepsilon$, then
\begin{eqnarray*}
\phi_{N-1}(\sqrt{N}x) &=& \frac{2^{1/4}}{\pi^{1/2}N^{{1}/{4}}}
\frac
{1}{(\sin\omega)^{{1}/{2}}} \sin
\biggl(\biggl(\frac{N}{2}+\frac
{1}{4}\biggr) (\sin2\omega- 2\omega)
+ \frac{3\pi}{4} \biggr) \\
&&{}\times\bigl(1 + O\bigl(N^{-1}\bigr)\bigr).
\end{eqnarray*}
\item[(4)]
If $x>\sqrt{2}+\delta$,
\[
\label{PRA5} \phi_{N-1}(\sqrt{N}x) = \frac{e^{-N I_1(x)}}{\sqrt
{4\pi\sqrt{2N}}}h(x) \bigl(1 +
O\bigl(N^{-1}\bigr)\bigr).
\]
\end{longlist}
\end{lemma}

\begin{pf*}{Proof of Lemma \ref{lemmanecessa}}
\textit{Part} (1): We can use the uniform asymptotics given by the exponential
region (4) in Lemma \ref{lemPR}. Precisely, by hypothesis, the
function $K(x):=ax^2+bx+I_1(x)$ is increasing in $[M,\infty)$, and by
Laplace's method,
\begin{eqnarray*}
I_N(M) &=& \int_{M}^{\infty}
\frac{e^{-N( ax^2+bx + I_1(x))}}{\sqrt{4\pi\sqrt{2N}}}h(x) \bigl(1
+ O\bigl(N^{-1}\bigr)\bigr) \,\d x
\\
&=& \frac{e^{-NK(M)}}{N|K'(M)|\sqrt{4\pi\sqrt{2N}}}h(M) (1+ O\bigl
(N^{-1}\bigr).
\end{eqnarray*}

\textit{Part} (2): Choose $\delta< \delta_0$ such that $-\sqrt{2}<M<\sqrt
{2}-\delta$. We eqnarray the integral $I_N(M)$ into three parts,
%
%
%e3.6 #&#
\begin{equation}
I_N(M) = \biggl(\int_{M}^{\sqrt{2}-\delta} + \int
_{\sqrt{2}-\delta
}^{\sqrt{2}+\delta} + \int_{\sqrt{2}+\delta}^{\infty}
\biggr):= I_1(M) + I_2 + I_3.
\end{equation}
We will show that the main contribution in this case comes from
$I_1(M)$. As in part~(1), it is easy to see that
%
%
%e3.7 #&#
\begin{equation}
I_3 = O\bigl(e^{-NK(\sqrt2)}\bigr).
\end{equation}
Next since $|x|^{1/4}|\operatorname{Ai}(x)|$ and
$|x|^{-1/4}|\operatorname{Ai}'(x)|$
are bounded functions on $\R$, a change of variables $z=I_1(-x)$ when
using part (2) of Lemma \ref{lemPR} immediately implies that for any
$\varepsilon>0$,
%
%
%e3.8 #&#
\begin{equation}
I_2 =O\bigl(e^{-N(a(\sqrt{2}-\delta)^2+b(\sqrt{2}-\delta))
+\varepsilon}\bigr).
\end{equation}

Now we estimate $I_1(M)$. Using the uniform asymptotics of $\phi_{N-1}$
we need to evaluate
%
%
%e3.9 #&#
\begin{eqnarray}
&&\frac{2^{1/4}}{\pi^{1/2}N^{{1}/{4}}} \int_{M}^{\sqrt{2}-\delta}
e^{-N( ax^2+bx)} \frac{1}{(\sin\omega)^{{1}/{2}}}
\nonumber
\\[-8pt]
\\[-8pt]
\nonumber
&&\hspace*{76pt}{}\times \sin\biggl
(\biggl(\frac
{N}{2}-
\frac{1}{4}\biggr) (\sin2\omega- 2\omega) + \frac{3\pi}{4}
\biggr)\,\d x.
\end{eqnarray}
Performing the change of variables $x = \sqrt{2} \cos\omega,
0<\omega
<\pi$ the integral above becomes (for some different $\delta>0$)
%
%
%e3.10 #&#
\begin{eqnarray}
\label{aaa333} &&\sqrt{2} \int_{\iota(M)}^{\pi-\delta}
e^{- N(2a\cos^2\omega+ \sqrt
{2}b\cos\omega)}
\nonumber
\\[-8pt]
\\[-8pt]
\nonumber
&&\hspace*{20pt}\qquad{}\times(\sin\omega)^{{1}/{2}} \sin\biggl(\biggl(
\frac
{N}{2}-\frac{1}{4}\biggr) (\sin2\omega- 2\omega) +
\frac{3\pi}{4} \biggr)\,\d\omega
\end{eqnarray}
for $\iota(M) = \arccos(2^{1/2}M)$.
We now rewrite $\cos^2\omega= \frac{1+\cos{2\omega}}{2}$ and use the
substitution $2\omega= z$ to obtain the integral
%
%
%e3.11 #&#
\begin{eqnarray}
\label{aaa3331} &&\frac{1}{\sqrt{2}} \int_{2\iota(M)}^{2\pi-2\delta}
e^{-N( a+a\cos{z} +
({b}/{\sqrt{2}})\cos({z}/{2}))}
\nonumber
\\[-8pt]
\\[-8pt]
\nonumber
&&\hspace*{31pt}\qquad{}\times \biggl(\sin\frac{z}{2}\biggr
)^{{1}/{2}} \sin
\biggl(\biggl(\frac{N}{2}-\frac{1}{4}\biggr) (\sin z - z) +
\frac{3\pi}{4} \biggr)\,\d z.
\end{eqnarray}
Last, we write
%
%
%e3.12 #&#
\begin{eqnarray}
\label{sinusexp} &&\sin\biggl(\biggl(\frac{N}{2}-\frac{1}{4}\biggr)
(\sin z - z) + \frac{3\pi}{4} \biggr)
\nonumber
\\[-8pt]
\\[-8pt]
\nonumber
&&\qquad= \frac{1}{2i}
\bigl[e^{i({N}/{2})(\sin z - z)} e^{if_1(z)} - e^{-i(
{N}/{2})(\sin z - z)} e^{-if_1(z)}
\bigr],
\end{eqnarray}
where $f_1(z) = -\frac{1}{4}(\sin z - z) + \frac{3\pi}{4}$.

Therefore, we just need to evaluate the asymptotics of
%
%
%e3.13 #&#
\begin{equation}
\int_{2\iota(M)}^{2\pi-2\delta} e^{- N m(z)}j(z) \,\d z,\qquad \int
_{2\iota(M)}^{2\pi-2\delta} e^{- N n(z)}k(z) \,\d z,
\end{equation}
where $m$ and $n$ are entire functions given by
%
%
%e3.14 #&#
%e3.15 #&#
\begin{eqnarray}
\label{eqm} m(z)&=&a+a\cos{z} + \frac{b}{\sqrt{2}}\cos\frac{z}{2} -
\frac
{i}{2}(\sin z - z),
\\
n(x)&=&a+a\cos{z} + \frac{b}{\sqrt{2}}\cos\frac{z}{2} +
\frac
{i}{2}(\sin z - z)
\end{eqnarray}
and $j(z) = \sin(\frac{z}{2})^{{1}/{2}}e^{if_1(z)}$, $k(z)=\sin
(\frac{z}{2})^{{1}/{2}}e^{-if_1(z)}$.

We will change our contour of integration and apply Laplace's integral
in the appropriate integrals. Notice that the steepest descent paths
are given by the equations
\begin{eqnarray*}
\label{cob} \Im\bigl(m(z)\bigr) &=& \sin x \biggl(a\sinh y + \frac
{\cosh y}{2}
\biggr) +\frac
{b}{\sqrt{2}}\sin\frac{x}{2}\sinh\frac{y}{2}-
\frac{x}{2} = \mbox{constant},
\\
\Im\bigl(n(z)\bigr) &=&\sin x \biggl(a\sinh y - \frac{\cosh y}{2}
\biggr) +
\frac
{b}{\sqrt{2}}\sin\frac{x}{2}\sinh\frac{y}{2}+
\frac{x}{2} = \mbox{constant}.
\end{eqnarray*}

The phase diagram for the steepest paths of $m$ is described as
follows. First all lines $x=2k \pi$, $k \in\N$ are steepest paths.
Second, for every $t\in(0,2\pi)$ the steepest path that passes through
$t$ goes from $0 - i \infty$ to $\pi+ i \infty$ if $b>0$ and from
$\pi- i \infty$ to $\pi+ i \infty$ if $b=0$. The real part of $m(z)$
is given by
\begin{eqnarray*}
\label{cab} \Re\bigl(m(z)\bigr) &=& \cos x \biggl(a\cosh y + \frac{1}{2}
\sinh y \biggr) + a +\frac{b}{\sqrt{2}}\cos\frac{x}{2}\cosh y -
\frac{y}{2},
\\
\Re\bigl(n(z)\bigr) &=& \cos x \biggl(a\cosh y - \frac{1}{2}\sinh y
\biggr) + a + \frac{b}{\sqrt{2}}\cos\frac{x}{2}\cosh y +
\frac{y}{2}.
\end{eqnarray*}

If we integrate $m(z)$ between two points $\alpha, \beta\in(0,2\pi
)$, we can deform our contour to be equal to the two steepest paths
that connect $\alpha$ and $\beta$ to $z= 0 - i\infty$. Precisely, we
deform our contour into three pieces: we first follow the steepest
descent path from $\alpha$ to a point with imaginary part $y_0<0$,
$|y_0|$ large. From there we go along the straight line $y=y_0$ until
we reach the steepest path that passes through $\beta$, $\gamma_{y_0}$,
and then we integrate on this steepest path back to $\beta$. We see
that if we choose $|y_0|$ large enough, every point in the straight
segment $y=y_0$ that we cross has real part $x$ sufficiently close to
$0$ so $\cos x >0$. This together with $ a>1/2$ implies that $\Re
(m(z))$ diverges to infinity as $y$ goes to negative infinity. The
trivial bound
%
%
%e3.16 #&#
\begin{equation}
\biggl|\int_{\gamma_{y_0}} e^{-Nm(z)} j(z) \,\d z\biggr| \leq\int
_{\gamma_{y_0}} e^{-N \Re(m(z))} \,\d z \sup_{z\in\gamma_{y_0}}
\bigl|j(z)\bigr|
\end{equation}
combined with the bounded length of $\gamma_{y_0}$ show that the
contribution of this part can be made as small as we want by choosing
$y_0$ large enough.

In the two remaining paths the imaginary part of $m$ is constant and
therefore we can apply Laplace's method to get the asymptotic behavior.
Since we assumed that $M < \sqrt{2}$ the contribution at $2\pi
-2\delta$
is negligible compared to the one at $2\iota(M)$. Indeed, by formula
(7.2.11) of \cite{Bleistein},
%
%
%e3.17 #&#
\begin{eqnarray}
\label{saast}&& \int_{2\iota(M)}^{2\pi-2\delta} e^{- N m(z)}j(z)
\,\d z
\nonumber
\\[-8pt]
\\[-8pt]
\nonumber
&&\qquad= \frac{e^{- N
m(2\iota(M)) + i (\pi- \alpha(M))}j(2\iota(M))}{N|m'(2\iota
(M))|}\bigl(1+O\bigl(N^{-1}\bigr)\bigr),
\end{eqnarray}
where $\alpha(M)$ is the angle of the steepest descent path of $m$ at
$z=2\iota(M)$,
%
%
%e3.18 #&#
\begin{equation}
\alpha(M) = \arctan\biggl(\frac{1-a \cos z}{2a\sin z + {b
\sin
(x/2)}/{\sqrt{2}}} \biggr).
\end{equation}
The above argument adapted to the function $n$ implies
%
%
%e3.19 #&#
\begin{eqnarray}
\label{saast2}&& \int_{2\iota(M)}^{2\pi-2\delta} e^{- N n(z)}k(z)
\,\d z
\nonumber
\\[-8pt]
\\[-8pt]
\nonumber
&&\qquad= \frac{e^{- N
n(2\iota(M)) + i (\pi- \alpha(M))} k(2\iota(M))}{N|n'(2\iota
(M))|}\bigl(1+O\bigl(N^{-1}\bigr)\bigr).
\end{eqnarray}

Noting that for any $x \in(0,2\pi)$ $|n'(x)| = |m'(x)|$, we can
combine \eqref{sinusexp}, \eqref{saast} and \eqref{saast2} to recover
that $I_1(M)$
is asymptoticly equivalent to
%
%
%e3.20 #&#
\begin{eqnarray}\qquad
&&\frac{2^{1/4}}{\pi^{1/2}N^{{5}/{4}}}\frac
{e^{-N(aM^2+bM)}}{2|m'(2\iota(M))|(\sin\omega)^{1/2}}
\nonumber
\\[-8pt]
\\[-8pt]
\nonumber
&&\qquad{}\times \sin\biggl[
\biggl(\frac{N}{2}-
\frac{1}{4} \biggr) (\sin2\omega- 2\omega) + \frac
{3\pi}{4} +
\alpha(M) \biggr]\bigl(1 + O\bigl(N^{-1}\bigr)\bigr).
\end{eqnarray}
This ends the proof of part (2) of lemma. The proof of part (3) follows
from the proof of part (2) and Laplace's method as in part (1) applied
to the integral
\[
\int_{M}^{-\sqrt{2}-\delta} e^{ax^2+bx+I_1(-x)} h(x) \,\d x = O
\bigl(e^{-N\lambda(M,a,b)}\bigr).
\]
In this case,
%
%
%e3.21 #&#
\begin{equation}
\label{eqclem82} L(M,a,b) = \frac{\sqrt{2\pi}h(\lambda
(M,a,b))}{2\lambda(M,a,b)+b +
I_1'(\lambda(M,a,b))}.
\end{equation}
We leave the details to the reader.
\end{pf*}

%In this case we have just as before:
%x^2} \d x.
%
%Now we use the fact that
%where $Z$ is an independent standard Gaussian random variable.
%So we can write \eqref{Laaa} as
%
%If we make $\alpha= 0$ and use \eqref{HFd} we recover the formula
%given in Theorem \ref{theo1euler}.
%
We now turn to the proof of Theorem \ref{meaneulerasym}.
\begin{pf*}{Proof of Theorem \ref{meaneulerasym}}
We can rewrite \eqref{eqmix} as
%
%
%e3.22 #&#
\begin{eqnarray}
\label{eqpartial} \E\chi(A_u) &= &(-1)^{N-1} \biggl(
\frac{\nu''}{\nu'} \biggr)^{
{(N-1)}/{2}} c(N,\nu)
\nonumber
\\
&&{}\times\int_{-\infty}^{\infty} \int_{-\infty
}^{{u}/{(\sqrt{2\nu''})}}
\phi_{N-1} \bigl(\sqrt{N}\bigl(\nu'x -\alpha y\bigr)
\bigr)\\
&&\hspace*{67pt}\quad{}\times e^{-N \nu''(x^2+y^2)} e ^{ ({N}/{2}) (\nu'x -\alpha y)^2}
\,\d x \,\d y,\nonumber
\end{eqnarray}
where
%
%
%e3.23 #&#
\begin{equation}
\label{cnnu} c(N, \nu) = 2 \nu''\bigl([N-1]!\sqrt{
\pi}\bigr)^{1/2} \frac{2^{-({(N-1)}/{2}})N}{\sqrt{\pi}\Gamma({N}/{2})}.
\end{equation}
For the case $\alpha\neq0$, we can change variables $z=\nu'x -\alpha
y$, $w = \alpha x + \nu'y$ to get
\[
x = \bigl(\nu'z + \alpha w \bigr) \biggl( \frac{1}{\alpha^2 + \nu'^2}
\biggr),\qquad  y = \bigl(\nu'w - \alpha z \bigr) \biggl(
\frac{1}{\alpha^2
+ \nu'^2} \biggr),
\]
and the above double integral becomes (using $\alpha^2 = \nu'' + \nu
' -
\nu'^2$)
% \frac{1}{\nu'' + \nu'}\int\int_{\nu'z + \alpha w \leq(\nu'' +
%which we rewrite as
%
\[
\frac{1}{\nu'' + \nu'}\int\int_{\nu'z + \alpha w \leq(\nu'' +
\nu
'){u}/{{2\nu''}}} \phi_{N-1} (
\sqrt{N}z ) e^{-{(N
\nu''(z^2+w^2))}/{(\nu'' + \nu')} } e^{N {z^2}/{2}} \,\d z \,\d w.
\]

So we have to evaluate the asymptotic behavior of the following integral:
\begin{eqnarray*}
J&=&\int_{-\infty}^{\infty} \phi_{N-1} ( \sqrt{N}z
) e^{{N (\nu
'-\nu'')z^2}/{(2(\nu'' + \nu'))} }\\
&&\hspace*{21pt}{}\times \int_{-\infty}^{{1}/{\alpha}
(
{(\nu''+\nu')u}/{(\sqrt{2\nu''})}- \nu' z)}
e^{-{N\nu''w^2}/{(\nu'+\nu
'')}} \,\d w \,\d z.
\end{eqnarray*}
We write the outside integral $ \int_{-\infty}^{\infty} \,\d z $ as
$\int_{-\infty}^{M} + \int_{M}^{\infty}$ with
\[
M = \frac{(\nu'+\nu'')u}{\sqrt{2\nu''}\nu'}.
\]
The inside integral is just a Gaussian integral, and therefore after a
straight-forward computation, the problem amounts to computing the
asymptotics of the two following one-dimensional integrals:
\begin{eqnarray*}
J_1 &=& \int_{M}^{\infty}
\phi_{N-1}(\sqrt{N}z) e^{-N({(\nu'^2+\nu
''-\nu')}/{(2(\nu''+\nu'-\nu'^2))} z^2 - {\sqrt{2\nu''}\nu
'u}/{(\nu''+\nu
'-\nu'^2)}z)} \,\d z,
\\
J_2&=& \int_{M}^{\infty}
\phi_{N-1}(\sqrt{N}z) e^{-N{(2(\nu'+\nu''))}/{(\nu''-\nu')}z^2}
\,\d z
\end{eqnarray*}
as $J = (J_1 + J_2) (1+O(N^{-1/2}))$ if $N$ is even and $J = (J_1 -
J_2) (1+O(N^{-1/2}))$ if $N$ is odd. Take $u \leq0$. We use Lemma \ref
{lemmanecessa} in both cases. Note that by \eqref{eq2p},
\[
a = \frac{\nu'^2+\nu''-\nu'}{2(\nu''+\nu'-\nu'^2)} > \frac
{1}{2} \quad\mbox{and}\quad b= - \frac{\sqrt{2\nu''}\nu'u}{\nu''+\nu'-\nu'^2}
\geq0.
\]

Now the condition $M\leq-\sqrt{2}$ ($M>-\sqrt{2}$) is exactly the
condition $u\leq-E_\infty'$ ($u>-E_\infty'$). Applying the appropriate
cases of Lemma \ref{lemmanecessa} we see that the integral $J_2$ is
negligible compared to $J_1$. A comparison with \eqref{eqqsela} and
\eqref{rafuta} gives the proof of part (1) and part (2) of the theorem
with $a$ and $b$ as above,
%
%
%e3.24 #&#
\begin{eqnarray}
\label{finalalpha} \alpha(w)&=& \arctan\biggl(\frac{1-a \cos\omega
}{2a\sin\omega+ {b
\sin(\omega/2)}/{\sqrt{2}}} \biggr),
\nonumber
\\[-8pt]
\\[-8pt]
\nonumber
 f(
\omega)&=& \bigl(\bigl|m'(2\omega)\bigr| \sin^{1/2} \omega
\bigr)^{-1}
\end{eqnarray}
and
%
%
%e3.25 #&#
\begin{equation}
\label{eqCt42} C(N,\nu,u) = \frac{1}{\nu''+\nu'} c(N,\nu) L(M,a,b),
\end{equation}
where $m$ is given in \eqref{eqm}, $c(N,\nu)$ in \eqref{cnnu} and
$L(M,a,b)$ in \eqref{eqclem82}.

If $\alpha= 0$, then the integral with respect to $y$ in \eqref
{eqpartial} can be explicitly computed and the mean Euler
characteristic is a single integral of the form given in Lemma~\ref
{lemmanecessa}. Applying part (1) and (2) of Lemma~\ref{lemmanecessa},
we get Theorem \ref{meaneulerasym} with
\[
C(N,\nu,u) = \frac{1}{\sqrt{2\pi}\Theta_\nu'(u)}h\biggl( \frac
{(\nu'+\nu
'')u}{\sqrt{2\nu''}\nu'}\biggr).
\]

Part (3) follows from symmetry of the Hamiltonian and \eqref{eqeulde}.
\end{pf*}

%s4 #&#
\section{Connection to mean field spin glasses}\label{secSP}

In this section we discuss our main motivation to study the problems
addressed in this manuscript. The function $H_N$ is the Hamiltonian of
a classical model in statistical physics, the mixed spherical $p$-spin
model \cite{CriSommers}. The study of the landscape of these
Hamiltonians is intimately related to the study of the most important
question in these systems, the $N$ limit of the Gibbs measure
\[
G_N (\bolds\sigma) = \frac{1}{Z_N} e^{\beta H_N(\bolds\sigma
)}.
\]

These mean-field models, as well as other spin glass models, are
well-known to be very challenging to analyze. It is believed (see \cite
{Leuzzi} and the references therein) that a subset of the spherical
models that we study here share the same interesting static and
dynamical behavior as the famous Sherrington--Kirkpatrick model at low
temperature. %

The understanding of the landscape of these Hamiltonians might prove
useful for the study of both static and dynamical questions of these
models. First, the structure derived from Theorem \ref{maintheorem} and
described below may shed a light on the metastability of Langevin
dynamics (in longer time scales than those studied in~\cite{BG97}).
Second, it may provide an insight (discussed below) in a possible
prediction for the structure of the Parisi measure, the functional
order parameter of these models.

The complexity of critical points $\theta_{k,\nu}(u)$ of finite index
has two pieces for negative values of $u$: one ``with a branching'' for
$u\in(-\infty, -E_\infty')$, another with a single curve, $u \in
(-E_\infty',0)$; see Figure~\ref{fig1}. This difference allows us to eqnarray
the models of Gaussian smooth functions on the sphere in two classes
that we describe now.

Let
%
%
%e4.1 #&#
\begin{equation}
G\bigl(\nu',\nu''\bigr):= \log
\frac{\nu''}{\nu'} -\frac{(\nu''-\nu') (\nu''-\nu'
+ \nu'^2)}{\nu'' \nu'^2} = \theta_{0,\nu}(-E_\infty).
\end{equation}

%
%de4.1 #&#
\begin{definition}
A mixture $\nu$ is called a \textit{pure-like} mixture if and only if
$G(\nu',\nu'') > 0$. If $G(\nu',\nu'')<0$, $\nu$ is called a
\textit
{full mixture}. When $G(\nu',\nu'')=0$, $\nu$ is called \textit{critical}.
\end{definition}

%
%ex1 #&#
\begin{example}
One can easily verify that all pure $p$-spins, $\nu(x)=x^p$, $p\geq3$
are pure-like while the spherical SK model, $p=2$, is critical.
\end{example}

\begin{figure}

\includegraphics{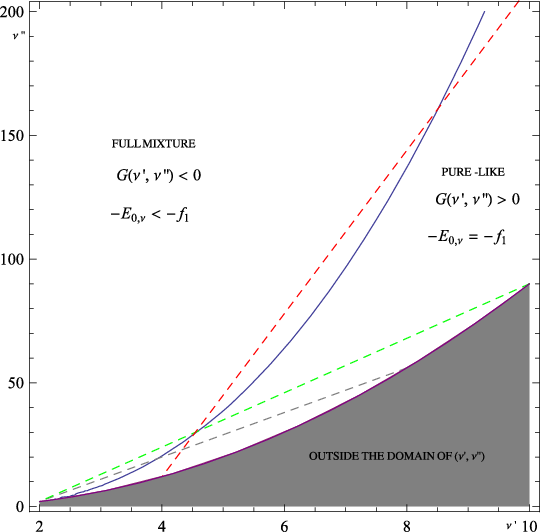}

\caption{Graph of $\nu' \times\nu''$. In blue, the level set $G(\nu
',\nu'')=0$, that is, the case where $\nu$ is critical. Dotted lines
are the possible values of $(\nu',\nu'')$ for the mixtures $2+6, 2+10$
and $4+30$. The gray region is outside the domain of possible values
for $(\nu',\nu'')$.}\label{zonas}% }
\end{figure}

%
%ex2 #&#
\begin{example}\label{example}
Consider the case
%
%
%e4.2 #&#
\begin{equation}
\label{poichi} \nu(t) = \mu t^2 + (1-\mu) t^p,
\end{equation}
where $\mu\in[0,1]$. Then, if $p > 3$, then it is possible to show
that there exists a $ 0<\mu_c(p)<1$ such that $\nu$ is \textit
{pure-like} if and only if $\mu\leq\mu_c(p)$.
$\mu_c(p)$ is given as the unique zero in $(0,1)$ of
\begin{eqnarray*}
&&-\frac{(p^2-2p)(1-\mu) (2(p^2-p)-3 (p^2-2p)\mu+(p-2)^2 \mu
^2 )}{2((p^2-p)(1-\mu)+2 \mu) (p+2\mu-p \mu)^2}
\\
&&\qquad{}+ \frac{1}{2}\log\biggl[1+p-\frac{2p}{p+2 \mu-p \mu} \biggr];
\end{eqnarray*}
see Figure~\ref{zonas}.
Remarkably, $p = 3$ in \eqref{poichi} is the only case where the
mixture is a pure-like mixture for all values of $\mu$.
\end{example}

%
%f2 #&#

It follows directly from the definition of pure-like and \eqref{defEk} that:
%
%
%pr5 #&#
\begin{proposition}\label{theoasd}
If $\nu$ is a pure-like mixture, then the sequence $E_k(\nu)$ is
strictly decreasing, and $E_k(\nu)$ converges to $E_{\infty}^{+}$ as
$k$ goes to infinity.
\end{proposition}

This proposition combined with Theorem \ref{tnofiniteindex} says if
the mixture $\nu$ is \textit{pure-like}, then the landscape of $\nu$ at
low levels of energy is similar to the pure case as in \cite{ABC}. In
particular, the same interesting layered structure for the
lowest critical values of the Hamiltonian $H_{N}$ holds. Namely, the
lowest critical values above the ground state energy are (with an
overwhelming probability) only local
minima, this being true up to the value $-NE_1(\nu)$, and that in a layer
above, $(-NE_1(\nu), -NE_2(\nu))$, one finds only critical values with
index 0
(local minima) or saddle points with index $1$, and above this layer one
finds only critical values with index $0, 1$ or $2$, etc.

There is one curiosity about pure-like mixtures. Define
%
%
%e4.3 #&#
\begin{equation}
\label{princvarpar} f_1:= \inf_{(a,b) \in[0,\infty)^2} \biggl\{
\frac{1}{2} \biggl(b + \nu'a + \frac{1}{b}\biggl(\log
\frac{a+b}{a}\biggr) \biggr) \biggr\}.
\end{equation}
%
%
%pr6 #&#
\begin{proposition}\label{tris11} $\nu$ is pure-like or critical if and
only if $f_1 = E_0(\nu)$.
\end{proposition}

The curiosity is that $f_1$ can be interpreted as the zero-temperature
limit of the 1-RSB Parisi functional in analogy to equation $(5.11)$ in
\cite{ABC}. We refer the reader to \cite{Talagrand} or Section~5 of
\cite{ABC} for a definition of this terminology. This leads us to the
following question:

%
%qu4.1 #&#
\begin{question} Is it true that a mixture is 1-RSB at low temperature
if and only if $\nu$ is pure-like?
\end{question}

The question raised above is consistent with a picture proposed by
physicists. In \cite{Leuzzi}, it is claimed that a $2+p$ spherical spin
glass model with $p\geq4$, at low temperature is either 1-RSB or its
Parisi measure has an absolute continuous part (a Full RSB or a 1-Full
RSB) depending on how much weight is assigned to the $2$-spin model.
The regions \textit{pure-like} and \textit{full mixture} seem to
numerically agree and to extend (since we do not need the $2$ spin
component) the one proposed by~\cite{Leuzzi}.

We end this section with the following statement about full-mixtures.
We first need the following result about the global minima of $H_N$
which is also of independent interest.

%
%th4.1 #&#
\begin{theorem}\label{Msri1} The following limit exists almost surely:
%
%
%e4.4 #&#
\begin{equation}
\lim_{N\rightarrow\infty} \frac{1}{N}\min_{\bolds\sigma}
H_N(\bolds\sigma):= -f_\infty.
\end{equation}
\end{theorem}
The following is now a corollary of Proposition \ref{tris11} and
Theorem \ref{maintheorem}.
%
%
%co4.1 #&#
\begin{corollary} If $\nu$ is a full mixture, then for any $u \in
(-E_0(\nu), -f_\infty)$, the probability of having a critical value
below $u$ goes to zero while the mean number of local minima is
exponentially large in $N$. Namely
for such $u$ there exist constants $0< C_1 < C_2$ such that for $N$
sufficiently large,
%
%
%e4.5 #&#
\begin{equation}
\E\Crt_{N,0}(-\infty, u) \geq e^{N C_1}\quad \mbox{and}\quad \Pro
\bigl(\Crt_N(-\infty,u)> e^{N C_1} \bigr) \leq
e^{-N C_2}.
\end{equation}
\end{corollary}

%s5 #&#
\section{\texorpdfstring{Proofs from Section \protect\ref{secSP}}
{Proofs from Section 4}}\label{sec5}
In this section we prove Propositions \ref{theoasd}, \ref{tris11} and Theorem
\ref
{Msri1}. We start by proving Theorem \ref{Msri1}. We will need to
introduce some notation and the lemma below. Let $\bolds\sigma^*$
be a point on the sphere such that $H_{N}(\bolds\sigma^*) = \min
_{\bolds\sigma} H_N$, and let $d$ denote the geodesic distance on
the sphere.
For $\rho, \alpha, K >0$, let
\[
B_{N,\rho} \equiv\bigl\{\bolds\sigma\in S_{N-1}(\sqrt{N})
\dvtx d\bigl(\bolds\sigma, \bolds\sigma^*\bigr) < \rho
\bigr\}
\]
and $A_{\varepsilon,\alpha,K}(N)$, be the event
%
%
%e5.1 #&#
\begin{equation}
A_{\varepsilon,\alpha,K}(N) \equiv\Bigl\{\sup_{\bolds\sigma\in
B_{N, \sqrt{N}\varepsilon}}
\bigl|H_{N} (\bolds\sigma) - H_{N}\bigl(\bolds\sigma^*
\bigr)\bigr| \leq K N \varepsilon^{\alpha} \Bigr\}.
\end{equation}

%
%le6 #&#
\begin{lemma}\label{smallball}
For any $0<\alpha<1$, there exist constants $K$, $K_1>0$ so that for all
$\varepsilon>0$ and all $N$ sufficiently large
%
%
%e5.2 #&#
\begin{equation}
\Pro\bigl(A_{\varepsilon,\alpha,K}(N)^c \bigr) < 2 e^{-K_1 N}.
\end{equation}
\end{lemma}
Note that this bound is independent of $\varepsilon$.
\begin{pf}
Clearly,
\[
A_{\varepsilon,K}(N) \supseteq\hat{A}_{\alpha,K}(N) \equiv\bigl\{ \|
H_{N} \|_{\alpha} \leq K N^{1-{\alpha}/{2}} \bigr\},
\]
where
%
%
%e5.3 #&#
\begin{equation}
\| H_{N} \|_{\alpha} = \sup_{\bolds\sigma, \bolds\sigma'}
\frac{|H_{N} (\bolds\sigma) - H_{N} (\bolds\sigma
')|}{d(\bolds\sigma,\bolds\sigma')^\alpha}.
\end{equation}

Now consider the centered Gaussian process $\mathbf{X}_{\alpha}$
field on $S_{N-1}(\sqrt{N})\times S_{N-1}(\sqrt{N})$ given by
%
%
%e5.4 #&#
\begin{equation}
\mathbf{X}_\alpha\bigl(\bolds\sigma, \bolds\sigma'\bigr) = %
\cases{ \displaystyle\frac{H_{N} (\bolds\sigma) - H_{N} (\bolds\sigma
')}{d(\bolds\sigma,\bolds\sigma')^\alpha},&\quad $\mbox{if } d
\bigl(\bolds\sigma,\bolds\sigma'\bigr)>0,$ \vspace*{2pt}
\cr
0, &\quad$\mbox{otherwise.}$} %
\end{equation}
Since the Gaussian field $H_{N} $ is $C^1$ almost surely, then
%
%
%e5.5 #&#
\begin{equation}
\label{unicaeq} \Pro\bigl(\hat{A}_{\alpha,K}(N)^c \bigr) = \Pro
\Bigl(\sup_{\bolds\sigma, \bolds\sigma'}\bigl |X_{\alpha
}\bigl(\bolds\sigma,
\bolds\sigma'\bigr)\bigr| > K N^{1-{\alpha}/{2}} \Bigr).
\end{equation}

But now a simple computation yields for %$\bolds\sigma_i \neq
$\bolds\sigma\neq\bolds\sigma'$,
%
%
%e5.6 #&#
\begin{eqnarray}
\E\mathbf{X}_{\alpha}^2\bigl(\bolds\sigma,
\bolds\sigma'\bigr)&=&\frac{2N}{d(\bolds\sigma
_1,\bolds\sigma'_1)^{2\alpha
}} \biggl[ 1-\nu\biggl(
\frac{1}{N}\bigl\< \bolds\sigma, \bolds\sigma'\bigr\>\biggr)
\biggr]
\nonumber
\\[-8pt]
\\[-8pt]
\nonumber
&=&\frac{2N}{(\sqrt{N}\theta)^{2\alpha}} \bigl[ 1-\nu(\cos
\theta) \bigr],
\end{eqnarray}
where $\theta$ is the angle between $\bolds\sigma, \bolds\sigma'$ in $\R^N$.

Therefore by the boundedness of $\nu'(x)$ in $[-1,1]$ there exists a
constant $C$ independent of $N$ such that [if $\alpha<1/2$ or $\alpha
<1$---using the boundedness of $\nu'(x)$ and $\nu''(x)$]
%
%
%e5.7 #&#
\begin{equation}
\sup_{(\bolds\sigma, \bolds\sigma')}\E\mathbf{X}_{\alpha
}^2\bigl(
\bolds\sigma, \bolds\sigma'\bigr) \leq CN^{1-\alpha}.
\end{equation}

Now, by Borell's inequality, (see pages 50 and 51 of \cite{AT07}, where
we take $u=K N^{1-{\alpha}/{2}}$, $\sigma_{T} \leq CN^{1-\alpha}$)
for all $\delta$, if $N, K$ is large enough
%
%
%e5.8 #&#
\begin{eqnarray}\quad
\Pro\Bigl(\sup_{\bolds\sigma, \bolds\sigma'}
\mathbf{X}_{\alpha}\bigl(
\bolds\sigma, \bolds\sigma'\bigr) > K N^{1-
{\alpha}/{2}}
\Bigr) &\leq& e^{\delta K N^{1-{\alpha}/{2}}}e^{
{-K^2 N^{2(1-{\alpha}/{2})}}/{(2CN^{1-\alpha})}}
\nonumber
\\[-8pt]
\\[-8pt]
\nonumber
&\leq& e^{-{K^2N}/{(4C)}}.
\end{eqnarray}
Taking $K_1= K^2/4C$ in the last equation, using \eqref{unicaeq} and
symmetry of $\mathbf{X}_{\alpha}$ the lemma is proven.
\end{pf}

\begin{pf*}{Proof of Theorem \ref{Msri1}}
Let $GS_N = \frac{1}{N}\min_{\bolds\sigma} H_N(\bolds\sigma
)$. We will show the existence of a constant $f_\infty$ so that for any
$\delta>0$ there exists $\varepsilon(\delta)$ such that if $N$ is large enough,
%
%
%e5.9 #&#
\begin{equation}
\label{caputio} \Pro\bigl( \bigl|GS_N +f_{\infty}\bigr|> \delta\bigr) \leq\Pro
\bigl(A_{\varepsilon(\delta),\alpha, K}(N)^c \bigr).
\end{equation}
The proof of Theorem \ref{Msri1} will then follow from \eqref{caputio}
and Borel--Cantelli's lemma since for all $\delta> 0$ by Lemma \ref
{smallball},
%
%
%e5.10 #&#
\begin{equation}
\sum_{N=1}^{\infty} \Pro\bigl( \bigl|GS_N
+f_{\infty}\bigr|> \delta\bigr) < \infty.
\end{equation}

We will prove \eqref{caputio} by showing that for any $\delta>0$ if $N$
is large enough $A_{\varepsilon,\alpha, K}(N) \subset\{ |GS_N
+f_{\infty
}|< \delta\}. $

On $A_{\varepsilon,\alpha,K}(N)$,
%
%
%e5.11 #&#
\begin{eqnarray}
\label{aburc} Z_{N,\nu}(\beta)&:=& \int_{S^{N-1}(\sqrt{N})}
e^{-\beta H_{N}(\bolds\sigma)} \Lambda_N(\d\bolds\sigma)
\nonumber
\\[-8pt]
\\[-8pt]
\nonumber
&\geq& e^{-\beta N GS_N - K \beta N \varepsilon^\alpha} \Lambda
_N(B_{N,\sqrt
{N}\varepsilon}).
\end{eqnarray}

Recall that $\Lambda_N(\d\bolds\sigma)$ is the surface measure
of $S_N(\sqrt{N})$ normalized to be a probability measure. We trivially
have the bound
%
%
%e5.12 #&#
\begin{equation}
\label{aburc2} \frac{1}{N} \log Z_{N,\nu}(\beta) \leq- \beta
GS_N.
\end{equation}

Combining \eqref{aburc} and \eqref{aburc2} we then have on
$A_{\varepsilon
,\alpha,K}(N)$,
%
%
%e5.13 #&#
\begin{eqnarray}
&&- \frac{1}{N \beta} \log Z_{N,\nu}(\beta) - K \varepsilon^{\alpha}
+ \frac
{1}{N \beta} \log\Lambda_N(B_{N,\sqrt{N}\varepsilon})
\nonumber
\\[-8pt]
\\[-8pt]
\nonumber
&&\qquad \leq
GS_N \leq-\frac{1}{N\beta} \log Z_{N,\nu}(\beta).
\end{eqnarray}

Note that using spherical coordinates and the inequality $\frac{2
\theta
}{\pi} \leq\sin\theta$ for $\theta\leq\frac{\pi}{2}$, we have for
$\varepsilon< \pi/2$,
%
%
%e5.14 #&#
\begin{eqnarray}
\Lambda_N(B_{N,\sqrt{N}\varepsilon}) &=& \biggl(\int_0^\varepsilon
\sin^{N-2} (\phi) \,\d\phi\biggr) \biggl(\int_0^\pi
\sin^{N-2} (\phi) \,\d\phi\biggr)^{-1}
\nonumber
\\[-8pt]
\\[-8pt]
\nonumber
&\geq&\biggl(\frac{2\varepsilon}{\pi}\biggr)^{N-1} \frac{1}{\pi(N-1)}.
\end{eqnarray}

So on $A_{\varepsilon,\alpha,K}(N)$, for some constant $C>0$,
%
%
%e5.15 #&#
\begin{equation}
- \frac{1}{N \beta} \log Z_{N,\nu}(\beta) - K \varepsilon^{\alpha}
+ C \varepsilon\leq GS_N \leq-\frac{1}{N\beta} \log
Z_{N,\nu}(\beta).
\end{equation}

By Holder's inequality the function $\frac{1}{N} \E\log Z_{N,\nu
}(\beta
)$ is convex in $\beta$, therefore its limit that we denote by
$F_\infty
(\beta)$ is also convex. The existence of this limit is given by the
famous Parisi formula \cite{Talagrand}, Theorem 1.1.

So $F(\beta)$ is convex, positive and grows at most linearly. This
easily implies that
%
%
%e5.16 #&#
\begin{equation}
\lim_{\beta\rightarrow\infty} \frac{1}{\beta} F_\infty(\beta) =
\sup_{\beta} \frac{1}{\beta} F_\infty(\beta):=
f_\infty\in[0,\infty).
\end{equation}

Therefore, for any $\delta_1 > 0$ one can take $N$ large enough so that
%
%
%e5.17 #&#
\begin{equation}
-\frac{F_{\infty}(\beta)}{\beta} - K \varepsilon^{\alpha} + C
\varepsilon-
\frac{\delta_1}{\beta} \leq GS_N \leq-\frac{F_\infty(\beta
)}{\beta} +
\frac{\delta_1}{\beta}.
\end{equation}
By taking $\beta$ large enough, part (a) of this theorem and by
choosing $\varepsilon$ sufficiently small, \eqref{caputio} is proven.
\end{pf*}

We now prove Propositions \ref{theoasd}.

\begin{pf*}{Proof of Proposition \ref{theoasd}}
If $\nu$ is pure-like, then $\theta_{k,\nu}(-E_\infty)>0$. Since
$\theta
_{k,\nu}(u)$ converges to negative infinity as $u$ goes to negative
infinity, $E_k(\nu)$ are well defined.
Furthermore, as $k$ goes to infinity, $\lambda_k^*(u)$ converges to
$-\sqrt{2}$ for any $u \leq- E_\infty$, implying that $\theta
_{k,\nu
}(u)$ converges to
$F(-\sqrt{2},u)$ pointwise.
Therefore, taking $u$ in a small neighborhood of $E_{\infty}^{+}$ and
using the fact that $\theta_{k,\nu}$ are increasing in that
neighborhood, we see that the zero of $\theta_{k,\nu}$ has to converge
to the zero of $F(-\sqrt{2},u)$. Namely $E_k(\nu)$ converges to
$E_{\infty}^{+}$.
\end{pf*}

%s5.1 #&#
\subsection{\texorpdfstring{Proof of Proposition \protect\ref{tris11}}
{Proof of Proposition 6}} We now
provide a proof
for Proposition \ref{tris11}. We will need a collection of calculus exercises.

%
%le7 #&#
\begin{lemma}\label{lemmadacon}
$f_1$ depends continuously on the first derivative $\nu'$.
\end{lemma}

%
%re5 #&#
\begin{remark}Note that while the $k$-complexity function depends
on the first two derivatives at $1$ of the covariance function $\nu$,
and $f_1$ depends only on the first derivative $\nu'$ and $E_0(\nu) =
f_1$ for any pure-like mixture.
\end{remark}
\begin{pf}
By solving for the critical points of \eqref{princvarpar}, we can get
an expression for $f_1$ in terms of $\nu'$. Namely,
\[
\label{f1function} f_1 = \frac{1}{2} \biggl(
\frac{\nu'y^2-1}{\nu'y} + \frac{1}{y} + \frac{\nu'y}{\nu
'y^2-1}\log\bigl(
\nu'y^2\bigr) \biggr) =y + \frac{\nu'-1}{y\nu'},
\]
where $y=y(\nu')$ is given by the unique solution of
\[
\label{eeqy} \biggl(\frac{\nu'y^2-1}{\nu'y} \biggr)^2y +
\frac{\nu'y^2-1}{\nu'y} = y\log\bigl(\nu'y^2\bigr), \qquad y>
\nu'^{-1/2}.
\]
In other words, $y=\frac{\sqrt{a}}{\sqrt{\nu'}}$ where $a$ is the
unique solution of
\[
\label{eeqa} a\log[a] - a + 1 - \frac{(a - 1)^2}{\nu'}=0,\qquad  a>1.
\]
This immediately implies the proof of Lemma \ref{lemmadacon}.
\end{pf}

%
%pr7 #&#
\begin{proposition}\label{calodoido} A mixture $\nu$ is critical if and
only if
%
%
%e5.18 #&#
\begin{equation}
f_1 = E_\infty= E_{0,\nu} = \frac{ \nu'' - \nu' + \nu'^2 }{\nu
'\sqrt{\nu''}}.
\end{equation}
\end{proposition}
\begin{pf}
If $\nu$ is critical, then $y=\frac{\sqrt{\nu''}}{\nu'}$ is the unique
solution of \eqref{eeqy} with $y > \frac{1}{\sqrt{\nu'}}$. Indeed,
\[
1-\frac{\nu''}{\nu'}+\frac{(-\nu'+\nu'') (-\nu'+\nu'^2+\nu''
)}{\nu'^3}-\frac{ (-1+{\nu''}/{\nu'} )^2}{\nu'} = 0.
\]
Plugging back the value of $y$ in \eqref{f1function} we get $f_1$.
On the other hand, if $f_1 = \frac{ \nu'' - \nu' + \nu'^2 }{\nu
'\sqrt{\nu''}}$, then one solves equation \eqref{f1function}\vspace*{-1pt} in $y$
to see
that the only positive solution is $y=\frac{\sqrt{\nu''}}{\nu}$. By the
definition of $y$ in \eqref{eeqy} this immediately implies that $\nu$
is critical. And trivially, $\nu$ critical is precisely the case where
$E_\infty= E_{0,\nu}$.
\end{pf}

Now we analyze the case where $\nu$ is critical or a full mixture, that
is, the case where $G(\nu',\nu'')\leq0$. In this case, the zero of the
complexity function can be explicitly computed and is given by
\[
-E_{0,\nu} = -E_{\infty}^{+},
\]
where $E_{\infty}^{+\bbook{}}$ was defined in \eqref{edocap1}. Note
that $E_{0,\nu}$ is a function of $\nu'$ and $\nu''$.

%
%pr8 #&#
\begin{proposition}\label{sabadao} If $G(\nu',\nu'')\leq0$, then
\[
\frac{\partial}{\partial\nu''} E_{0,\nu} = 0\quad \mbox{if and only
if}\quad G\bigl(
\nu',\nu''\bigr)=0.
\]
\end{proposition}

\begin{pf} Let
\[
A\bigl(\nu',\nu''\bigr)=\sqrt{\bigl(
\nu''-\nu'^2+
\nu'\bigr) \biggl(\bigl(\nu'+\nu''
\bigr) \log\biggl[\frac{\nu''}{\nu'} \biggr]-2 \bigl(\nu''-
\nu'\bigr) \biggr)}.
\]
Calculating the derivative $\frac{\partial}{\partial\nu''} E_{0,\nu}$
one gets
%
%
%e5.19 #&#
\begin{eqnarray}
\label{shabu}  &&\biggl(\nu'^2 \nu''
\bigl(\nu'+\nu''\bigr) \log\biggl[
\frac{\nu''}{\nu'} \biggr]+\bigl(\nu''-
\nu'\bigr)\biggr) \nonumber\\
&&\qquad{}\times\bigl(\nu'^3+
\nu''^2-\nu'^2
\bigl(1+3 \nu''\bigr)-2 \nu' \sqrt{\nu
''} A\bigl(\nu',\nu''
\bigr) \bigr)\\
&&\qquad{}\times \bigl(2 \nu'' \bigl(\nu'+
\nu''\bigr)^2 A\bigl(\nu',
\nu''\bigr) \bigr)^{-1}.\nonumber
\end{eqnarray}
Sufficiency comes from a simplification of the above formula. To get
necessity we solve a second degree equation on the variable $M = \log
[\frac{\nu''}{\nu'} ]$ to see that this second degree
equation has a unique zero given by
\[
\frac{\nu'^2 - \nu'^3 - 2 \nu' \nu'' + \nu'^2 \nu'' + \nu
''^2}{\nu'^2
\nu''}.
\]
This is precisely $G(\nu',\nu'')=0$.
\end{pf}

With the above propositions we now prove Proposition \ref{tris11}.

\begin{pf*}{Proof of Proposition \ref{tris11}}
If $\nu$ is critical, Proposition \ref{tris11} is Proposition \ref
{calodoido}. Now suppose that $\nu$ is pure-like. By Lemma \ref
{lemmadacon} and \eqref{complexityfunction2}, both $f_1(\nu):=f_1$ and
$E_0(\nu)$ are independent of $\nu''$. Consider then another mixture
$\mu$ such that $\mu'=\nu'$ and $\mu$ satisfies $G(\mu',\mu'')=0$.
Since $G$ is continuous on its domain, we have
\[
f_1(\nu)=f_1(\mu)=E_{0}(
\mu)=E_0(\nu).
\]
On the other hand, if $\nu$ is a full-mixture, Proposition \ref
{sabadao} combined with Lemma \ref{lemmadacon} shows that $f_1\neq
E_0(\nu)$. This ends the proof of Proposition \ref{tris11}.
\end{pf*}

\section*{Acknowledgments}
We want to underline our debt to Michel Ledoux for his friendly help
with the results of Section~\ref{sec5}. We also would like to thank
Jiri Cerny for a careful reading of this manuscript and Yan Fyodorov
for pointing out that the method used in this paper is similar to \cite
{fyodorov-2004-92} and \cite{Fyodorov}. We want to thank MSRI, IMPA,
Universit\'{e} de Marseille for their hospitality and the Universit\'
{e} de Nice where a mini-course based on these results were given. We
are also in debt to an anonymous referee who helped us to improve the
presentation of this manuscript.

% imsref loaded by akundreckaite, 2013-10-03 09:52:13
%

%
% zodis "Acknowledgments" paliekamas pagal autoriu

%suskaldyti doi

\printaddresses

\end{document}